\documentclass[12pt]{article}

\makeatletter
\renewcommand{\@biblabel}[1]{#1.}
\makeatother
\newtheorem{De_d}{\rm \bf Definition}
\newtheorem{Ex_d}{\rm \bf Example}
\newtheorem{Rem_d}{\rm \bf Remark}
\newtheorem{Sa_d}{\rm \bf Proposition}
\newtheorem{Le_d}[Sa_d]{\rm \bf Lemma}
\newtheorem{Th_d}[Sa_d]{\rm \bf Theorem}
\newtheorem{Fo_d}[Sa_d]{\rm \bf Corollary}
\newenvironment{De}{\begin{De_d}{$\!\!\!$\bf .}}{\end{De_d}}
\newenvironment{Ex}{\begin{Ex_d}{$\!\!\!$\bf .}}{\end{Ex_d}}
\newenvironment{Rem}{\begin{Rem_d}{$\!\!\!$\bf .}}{\end{Rem_d}}
\newenvironment{Sa}{\begin{Sa_d}{$\!\!\!$\bf .}}{\end{Sa_d}}
\newenvironment{Le}{\begin{Le_d}{$\!\!\!$\bf .}}{\end{Le_d}}
\newenvironment{Th}{\begin{Th_d}{$\!\!\!$\bf .}}{\end{Th_d}}
\newenvironment{Fo}{\begin{Fo_d}{$\!\!\!$\bf .}}{\end{Fo_d}}

\newcommand{\RR}{{\rm I\kern-0.14em R }}
\newcommand{\NN}{{\rm l\kern-0.14em N}}
\newcommand{\CC}{{\rm\raise 0.192ex\vbox{\hrule height
                  1.22ex width 0.8pt}\kern-0.29em C}}
\newcommand{\ZZ}{ {\sf Z}\hspace{-0.4em}{\sf Z}\ }

\begin{document}
\begin{flushleft}
{ \Large   Abel's Theorem
in the Noncommutative Case }
\end{flushleft}

\vspace{0.5cm}

\begin{flushleft}
{FRANK LEITENBERGER } \\
{\small \it Fachbereich Mathematik, Universit\"at Rostock,
Rostock, D-18051, Germany. \\
e-mail: frank.leitenberger@mathematik.uni-rostock.de }
\end{flushleft}


\vspace{1.0cm}

\noindent
short title: Abel's Theorem


\vspace{0.5cm}

\begin{flushleft}
\small {\bf Abstract.} We define noncommutative binary forms.
Using the typical representation of Hermite we prove the
fundamental theorem of algebra and we derive
a noncommutative Cardano formula for
cubic forms. We define quantized elliptic and hyperelliptic
differentials of the first kind. Following Abel
we prove Abel's Theorem.
\end{flushleft}

\setcounter{equation}{0}

\begin{flushleft}
\bf 1. Introduction
\end{flushleft}
Finding quantum deformations of the basic objects of algebraic
geometry is an open problem. There are different approaches
(cf. \cite{Kon,Man}). Plane curves with
genus $g\geq 1$ are related to
elliptic and Abelian integrals. Addition theorems of
these functions are the content of Abel's theorem. Our approach is
based on classical invariant theory and Abel's ideas.
We demonstrate that a $h$-deformation has no influence on the
validity of Abel's theorem for hyperelliptic differentials of
the first kind. For this purpose we develop a theory of a
$h$-deformation of the classical invariant theory of binary forms.

Let $C$ be a hyperelliptic curve of genus $g$ with
Weierstrass points $a_1,\ \cdots\  ,a_{2g+2}$, i.e. $C=\{(t,u)\in
\CC^2| f(t,u)=0   \}\cup \{ \infty \}  $, where
$f(t,u)=u^2-p(t)=u^2- (t-a_{1})\ ..... \ (t-a_{2g+2})$.
We consider $k$ points $(t_1,u_1),...,(t_k,u_k)$ on $C$
and form the sum of complex integrals
\[ S=
\int_{\infty}^{t_1,u_1} \frac{v(t)dt}{u}+ \int_{\infty}^{t_2,u_2}
\frac{v(t)dt}{u}+\ ...\ + \int_{\infty}^{t_k,u_k}
\frac{v(t)dt}{u},
\]
where $v(t)$ is a polynomial of degree $\leq g-1$ and a
complex path from $\infty$ to $(t_i,u_i)$ is chosen in some way.
We call $\frac{v(t)dt}{u}$ a hyperelliptic differential of
the first kind. We consider on $C$ the meromorphic function
\[  q(t,u,c_1,...,c_{m_1},d_1,...,d_{m_2} )=
a(t,c_1,...,c_{m_1})u-b(t,d_1,...,d_{m_2})  \] with polynomials
$a(t),b(t)$ depending rationally on parameters $c_i,d_i$. Abel's
theorem tells us that $S=0=const$  if the
$(t_1,u_1),...,(t_k,u_k)$ are the points of intersection of
$f(t,u)=0$ with the variable algebraic curve $q(t,u,c_i,d_i)=0$.
We remark that due to $u^2-p(t) \equiv 0$ on $C\backslash
\{\infty\}$ we can reduce the intersections of any variable
algebraic curve $\tilde{q}(t,u,\tilde{c_1},...,\tilde{c_m})=0$
depending rationally on the $\tilde{c}_i$ to the above situation.

Abel's theorem has the differential form
\begin{equation} \label{x0}
\epsilon_1 \frac{v(t_1)dt_1}{\sqrt{p(t_1)}}+ \epsilon_2
\frac{v(t_2)dt_2}{\sqrt{p(t_2)}}+ ...+ \epsilon_k
\frac{v(t_k)dt_k}{\sqrt{p(t_k)}}=0
\end{equation}
with $\epsilon_i = \pm 1$.

Now we replace the $t_i$ in equation (1) by the quotients of
homogeneous coordinates $\frac{x_i}{y_i}$. Furthermore we multiply
all components of equation (1) by a certain common factor (cf.
Remark 20 below). We can express the hyperelliptic differentials,
which are now endowed by this common factor, by bracket symbols
\[
(ij):=\left| \begin{array}{cc} x_i & x_j \\ y_i & y_j
\end{array} \right|,
\ \ \ \ \ \ \ \ \ \ (i,dj):=\left| \begin{array}{cc} x_i & dx_j \\
y_i & dy_j
\end{array} \right|
\]
introduced by Clebsch. Now the components are
$SL(2)$-invariant with respect to the natural action on the
$x$-$y$-plane. Therefore we can consider the differential form of
Abel's theorem as a proposition about the sum of certain
differential invariants. Clebsch introduced the invariant
theoretical descriptions of elliptic integrals in \cite{Cl},
p. 228. Further developments are contained in \cite{CG}.

Most of the modern treatments about Abel's theorem use the
Riemanian ideas. For our noncommutative generalization of the
differential form (\ref{x0}) we will follow the original proof of
Abel (cf. \cite{Ab}). Abel reduced equation (\ref{x0}) to an
algebraic identity, which is a consequence of an expansion into
partial fractions (cf. \cite{We}, p.28).

We introduce an algebra $H_I$ of noncommutative homogeneous
coordinates of points of a noncommutative line. The quantum
group $U_h(sl(2))$ plays the role of the Cayley-Aronhold
differential operators of the classical invariant theory
which generate a Lie Algebra $sl(2)$.
Because of the Ore property of $H_I$ we can extend
this algebra to a division algebra.

We define $h$-deformed bracket symbols by
$(ij)=x_iy_j-y_ix_j-hy_iy_j$. The classical invariant theory is
essentially determined by the algebra of symbols $(ij)$. We make
the observation that the algebra of symbols is isomorphic to the
classical case. Furthermore the two fundamental theorems of
invariant theory are valid.

We define $n$-forms and invariants of $n$-forms. It turns out that
the Clebsch Gordan symbolic method works in our situation. By
computer calculations we derive the simplest invariants and
covariants of quadratic and cubic forms. We use polarization
operators in order to replace the independent variables of
$n$-forms by the coordinates of arbitrary noncommutative points.

Using the theory of the typical representation of Hermite (cf.
\cite{He}) we derive the fundamental theorem of algebra for
noncommutative $n$-forms, i.e. we decompose $n$-forms into $n$
commuting linear forms in a certain skew field extension. We
describe explizite decompositions for $n=2,3$.

In order to prove Abel's theorem we combine the ideas of Abel and
Clebsch. We define elliptic and hyperelliptic differentials as
differential invariants in certain differential modules. The
addition theorem for elliptic differentials and Abel's theorem for
hyperelliptic differentials appear as identities in these
differential modules. The classical set of intersection points
$(t_i,u_i)$ is replaced by a set $(X_i,Y_i)$ of homogeneous
coordinates of the zeros of a noncommutative $k$-form $r$, which
corresponds to an elimination expression $r(t)$ of $f(t,u)$ and
$q(t,u,c_i,d_i)$ in the classical case. The analogy of the proof
to the classical case is based on the fact that the elements
$(ij),y_k,dz_l$ commute. The $h$-deformation has no influence over
any step of the classical invariant theoretical proof.

We used the computer algebra program Mathematica 3.0 to form the
computations.

\begin{flushleft}
\bf \large 2. Preliminaries
\end{flushleft}

In this section we introduce some basic concepts about the quantum
group $U_h(sl(2))$ and the braided module algebras $H_I$. We
follow the similar considerations in \cite{Le} for the quantum
group $U_q(gl(2))$.

Let $U_h(sl(2))$
with $ h\in  \CC $, $h\neq 0$
be the complex unital Hopf algebra,
determined by generators $E,F,H$ and relations
\begin{eqnarray} \label{eq1}
\begin{array}{lll}
 {[} H , E {]} & = &  E(\cosh (hF))+(\cosh (hF))E  \\
 {[} H , F {]} & = & -\frac{2}{h}  \sinh (hF)      \\
 {[} E , F {]} & = &  H.
\end{array}
\end{eqnarray}
The Hopf multiplication is given by
\begin{eqnarray*}
\begin{array}{ccccccccc}
\Delta (E) & = & E & \otimes  & e^{-hF} & + & e^{hF} & \otimes & E, \\
\Delta (F) & = & F & \otimes  & 1       & + &      1 & \otimes & F, \\
\Delta (H) & = & H & \otimes  & e^{-hF} & + & e^{hF} & \otimes & H.
\end{array}
\end{eqnarray*}
$U_h(sl(2))$ has the counit
\begin{eqnarray*}
\epsilon (E)=\epsilon (F)=\epsilon (H)=0, \ \ \ \ \
\epsilon (1)=1.
\end{eqnarray*}
and the  coinverse
\begin{eqnarray*}
S(E) = - e^{hF} E e^{-hF} ,\ \ \ \ \
S(F) = -        F         ,\ \ \ \ \
S(H) = - e^{hF} H e^{-hF}.
\end{eqnarray*}

For an arbitrary ordered index set $I$ we consider the complex
unital $U_h(sl(2))$-module algebra $A_I$, which is freely generated
by the variables $x_i$, $y_i$, $i \in I$. To determine the action
of $U_h(sl(2))$ on $A_I$, we set
\begin{eqnarray} \label{eq2}
\begin{array}{lll}
H 1   =  0                   ,\ \ \ \ \   &
H x_i =  x_i,            &
H y_i =  - y_i,\\
E 1   =  0                   ,\ \ \ \ \   &
E x_i =  0                   ,            &
E y_i =  x_i  \\
F 1   =  0                   ,\ \ \ \ \   &
F x_i =  y_i,            &
F y_i =  0                   ,
\end{array}
\end{eqnarray}
and require
\begin{eqnarray} \label{eq3}
\begin{array}{ccll}
E(ab) & = & E(a)  e^{-hF} (b) & + e^{hF}(a)  E(b), \\
F(ab) & = & F(a)  b           & +     a      F(b), \\
H(ab) & = & H(a)  e^{-hF}(b) & + e^{hF}(a)  H(b)
\end{array}
\end{eqnarray}
for $a,b\in H_I$. A proof can be given along the lines of
\cite{Kl}, p.19 where a similar module occurs.

Furthermore we consider the ideal $J$ of $A_I$, which is generated
by the elements
\begin{eqnarray*}
\begin{array}{lll}
x_iy_i - y_ix_i-hy_i^2, \\
x_jy_i - y_ix_j-hy_iy_j, \\
y_jx_i - x_iy_j+hy_iy_j, \\
y_jy_i - y_iy_j,   \\
x_jx_i - x_ix_j+hx_iy_j-hy_ix_j-h^2y_iy_j
\end{array}
\end{eqnarray*}
for $i<j$. One can show that the ideal $J$ is $U_h(sl(2
))$-invariant (i.e. $E(J),F(J),H(J)\subseteq J$
(cf. \cite{Ma}). Therefore the action of $U_h(sl(2 ))$ on $A_I$
induces an action on $H_I:=A_I/J$. In the following we identify
$x_i$ with its image under the quotient map $ A_I \rightarrow H_I
$. I.e. $H_I$ is the unital $U_h(sl(2 ))$-module algebra with
generators $x_i,y_i$ $i\in I$ and relations
\begin{eqnarray} \label{eq4}
\begin{array}{ll}
x_iy_i  =  y_ix_i+hy_i^2, & \\
x_jy_i  =  y_ix_j+hy_iy_j,\ \ \ \ \ \ \ \ \ \ \ \ &
y_jx_i  =  x_iy_j-hy_iy_j, \\
y_jy_i  =  y_iy_j, &
x_jx_i  =  x_i x_j-h x_i y_j+h y_i x_j+h^2 y_i y_j
\end{array}
\end{eqnarray}
for $i<j$.

\begin{Rem} \label{Rem1}  \rm
$H_I$ is a {\it braided} module algebra (cf. \cite{Ma}).
We have $(x_i+x_j)(y_i+y_j)=(y_i+y_j)(x_i+x_j)+h(y_i+y_j)(y_i+y_j)$.
The algebra $H_I$ has the
PBW property for an arbitrary order of the variables
$x_i,y_i$, $i\in I$.
We can consider the equations (\ref{eq4}) as rules
in order to express an element $a\in H_I$ by the PBW basis
                 $x_{i_1}^{n_1} y_{i_1}^{m_1}
                  x_{i_2}^{n_2} y_{i_2}^{m_2}
                         ...
                  x_{i_k}^{n_k} y_{i_k}^{m_k}$
with $i_1<i_2<...<i_k$ and $k\in \NN$.
The rules (\ref{eq4}) do not change the type of
homogeneity with respect to the indices $i\in I$.

We note that $H_I$ is
the only module algebra with quadratic relations
which is braided and has the PBW property.
\end{Rem}

We will extend $H_I$ to a division algebra $Q_I$.
Let $R$ be an algebra, $\sigma$ an endomorphism of $R$
and $\delta$ a $\sigma$-derivation with
\[ \delta (ab) = \delta (a) b+\sigma (a) \delta (b).\]
We say that $S=R[x;\sigma, \delta ]$
is an {\it Ore extension of $R$}
if $S$ is freely generated over $R$ by an element $x$
subject only to the relation
\[     xa = \sigma (a) x  +  \delta (a)   \]
for $a\in R$ (cf. \cite{Mc}).

\begin{Le} \label{Sa1}
The algebra $H_I$ is an iterated Ore extension.
\end{Le}

{\it Proof.}
We restrict the consideration to the finite generated algebras $H_I$
with $I=\{0,1,...,n \}$.
In the case of an infinite index set $I$ one can proceed by
transfinite induction.
We consider the tower of subalgebras
\[   \CC \subset H_0' \subset H_0   \subset H_1' \subset
           \ \cdots \
            \subset H_n' \subset H_n   \]
where $H_j' \cong H_{j-1} [y_{j},\sigma,\delta]$
is the subalgebra generated by $H_{j-1}$ and $y_j$,
where $\delta = 0$ and $\sigma$ is determined by
\[  \sigma (x_i)=   x_i-hy_i,\ \ \ \ \
    \sigma (y_i)=   y_i \]
for $i<j$.
Furthermore we have $H_j \cong H_j' [x_{j},\sigma,\delta]$
where $\sigma$ and $\delta$ are determined by
\[   \sigma (x_i)    = x_i + h y_i,\ \ \ \ \
     \sigma (y_i)    = y_i,\ \ \ \ \
     \sigma (y_{j})  = y_j,      \]
\[   \delta (x_i)    = -h(x_i-hy_i)y_j,\ \ \ \ \
     \delta (y_i)    =  hy_iy_j,\ \ \ \ \
     \delta (y_{j})  =  hy_j^2 \]
for $i<j$.
The verification of the isomorphies
is similar to \cite{Ka}, p. 81.
The assertion follows.
${\bullet}$

\begin{Sa} \label{Sa2}
The algebra $H_I$ has no zero divisors
and $H_I$
has an uniquely determined right quotient division ring $Q_I$,
whose elements are the right quotients
$ab^{-1}$ for $a,b\in R$, $b\neq 0$.
\end{Sa}

{\it Proof.}
Ore proved that a noncommutative ring without zero divisors
has an uniquely determined quotient division ring if
\begin{eqnarray} \label{eq5}
  aR \cap bR \neq \{ 0\}\ \ \ \ \ \ \ \forall a,b \in R,
                            \ \ \ \ \ a,b\neq 0,
\end{eqnarray}
i.e., two nonvanishing elements $a,b$
have a common right multiple (cf. \cite{Or}).
Curtis proved that an iterated Ore extension
of a ring without zero divisors with the property (\ref{eq5})
is again a ring without zero divisors with the property (\ref{eq5})
(cf. \cite{Cu}).
The assertion follows from both results together with Lemma \ref{Sa1}.
${\bullet}$

\begin{Ex} \label{Ex1} \rm
For $x_1$ and $x_2$ we find the common right multiple
\[    x_2 ( x_1-h y_1 )
    = x_1 ( x_2-h y_2 ).  \]
Therefore we can represent the left quotient
$x_1^{-1} x_2$
by the right quotient \\
$  ( x_2-h y_2 )( x_1-h y_1 )^{-1}   $.
\end{Ex}

\begin{Rem} \label{Rem2} \rm
We can consider the variables $x_i, y_i$
as noncommutative homogeneous coordinates
of different points of the complex projective line.
The noncommutative coordinate is given by
\[  z_i: =  x_i y_i^{-1}+\frac{h}{2}     \in Q_I. \]
For an involution $x_i^*:=x_i$, $y_i^*:=y_i$ and ${\rm Re}\ h=0$
we have $z_i^*=z_i$. This motivates the $\frac{h}{2}$\ -term.
For two points we have the commutation relation
\begin{eqnarray*}
z_jz_i=z_iz_j+2h(z_j-z_i), \ \ \ \ \ \forall i,j
\end{eqnarray*}
or $(z_j+h)(z_i-h)=(z_i+h)(z_j-h)$.
\end{Rem}

\begin{flushleft}
{\bf \large 3. The Invariants of the algebra $H_I$. }
\end{flushleft}

\begin{De} \label{De1}
{\rm We say that the element $v\in H_I$ is an
{\it invariant element}, if
\begin{eqnarray*}
E(v)=F(v)=H(v)=0.
\end{eqnarray*}
}\end{De}

\begin{Rem} \label{Rem3} {\rm
Alternatively, one can define invariant elements
by the Quantum group $SL_h(2  )$.
The algebra structure is given by
the generators $a,b,c,d$ and the relations \\
$ab-ba=h(1-a^2), \ \
ac-ca=hc^2,\ \
bd-db=h(d^2-1)$, \\
$cd-dc=-hc^2, \ \
ad-da=hc(d-a),\ \
bc-cb=h(dc+ca)$ \\
and $ad-bc+hac=1$.
A coaction $\omega$
of $SL_h(2 )$ on $H_I$ is given by
\begin{eqnarray*}
\left(
\begin{array}{c}
\omega (x) \\
\omega (y)
\end{array}
\right) :=
\left(
\begin{array}{c}
a \otimes x + b \otimes y \\
c \otimes x + d \otimes y
\end{array}
\right).
\end{eqnarray*}
Now we can define invariant elements by the equation
\begin{eqnarray*}
\omega (v) = 1 \otimes v.
\end{eqnarray*}
}\end {Rem}

By the rules $(2)$ sums and products of invariant
elements are again invariant elements.
Therefore they form a subalgebra
$H_I^{Inv} \subset H_I$.

\begin{Ex} \label{Ex2}\rm
The simplest invariant
element of $H_I$ is given by
\begin{eqnarray*}
(ij):=x_iy_j-y_ix_j-hy_iy_j.
\end{eqnarray*}
We have also
$(ij) = x_i y_j - x_j y_i =  y_j x_i - y_i x_j$.
\end{Ex}
In the following we will see,
that the subalgebra $H_I^{Inv}$
has many classical properties.
\begin{Sa} \label{Sa3}
The bracket symbols
$(ij)$, $i,j\in I$ form a central, commutative
subalgebra $S_I \subseteq H_I^{Inv}$.
\end{Sa}

\noindent
{\it Proof:}
The identities
$x_k (ij)  = (ij) x_k$,
$y_k (ij)  = (ij) y_k$,  $i,j,k,\in I$
can be checked by an explizite calculation.
Therefore the bracket $(ij)$ is central.
In particular we have
$(kl)(ij) = (ij)(kl)$, $\forall i,j,k,l \in I$.
Therefore $S_I$ is commutative. $\bullet$

The bracket symbols $(ij)$ have the following simple properties
$(ii) = 0$ and
\begin{eqnarray} \label{eq6}
(ji)+(ij)=0,\ \ \ \ \ \ \ \ (ij)(kl) + (ik)(lj) + (il)(jk)=0
\ \ \ \ \ \ \ \forall i,j,k,l.
\end{eqnarray}
The second identity is the Grassmann Pl\"ucker relation.

Now we consider the two Fundamental Theorems of invariant theory.
We have the First Fundamental Theorem:
\begin{Th} \label{Sa4}
We have $H_I^{Inv}=S_I$,
i.e. the algebra $H_I^{Inv}$
is generated by the bracket symbols $(ij)$.
\end{Th}
A proof can be given similar to the case of the $U_q
(gl_2)$-symmetry (cf. \cite{Le}).

Because $H_I^{Inv}$ is generated by bracket symbols,
we can consider $H_I^{Inv}$ as a quotient of an algebra $B_I$,
freely generated by commuting symbols $[ij]$
and an ideal $S$ of syzygies
\\ (i.e. $(ij)\cong [ij]+S$). The Second Fundamental Theorem
of classical invariant theory tells us that $S$ is generated
by the elements $[ij]+[ji]$ and $[ij][kl]+[ik][lj]+[il][jk]$
(cf. \cite{Ol}).

The Second Fundamental Theorem is also valid for $h\neq 0$:
\begin{Th} \label{Sa5}
Using only arithmetical operations of bracket symbols,
every vanishing expression of bracket symbols
can be  transformed into
\[ \sum_{ij}  ((ij)+(ji))I_{ij}
+\sum_{ijkl}  ((ij)(kl)+(ik)(lj)+(il)(jk))I_{ijkl} \]
with certain expressions $I_{ij},I_{ijkl}$ of bracket symbols.
\end{Th}

\noindent
{\it Proof:}
For every $h\in \CC$ we
consider the homomorphism
$g_h:B_I\rightarrow H_I^{Inv}=H_I^{Inv}(h)$
with $g_h([ij])=(ij)$.
Let $S_h:=Ker\ g_h$.
By the classical Second Fundamental Theorem
and (\ref{eq6})
we have $S_0 \subseteq S_h$.
Conversely, let $e\in B_I$.
We consider the PBW expansion of $g_h(e)$.
If we apply the rules (\ref{eq4}) to a monomial of the $x_i, y_i$
then the $x_i, y_i$ permute
without changing the coefficient of the monomial
and we obtain additional
terms with a higher degree in the $y_i$.
Therefore we have the decomposition
\[  g_h(e)=e_0^{(h)}+e_1^{(h)}  \]
where $e_0^{(h)}$  contains all monomials
$x_{1}^{n-i_1} y_{1}^{i_1}
 x_{2}^{n-i_2} y_{2}^{i_2}
       ...
 x_{n}^{n-i_n} y_{n}^{i_k}$
with $\sum_j i_j =\frac{1}{2} \sum_j n_j$
and $e_1^{(h)}$ contains all monomials with
$\sum_j i_j >\frac{1}{2} \sum_j n_j$.
The coefficients of $e_0^{(h)}$
are independent of $h$ and
$e_1^{(h)}$ is $0$ for $h=0$.
Therefore $g_h(e)=0$ yields $g_0(e)=0$,
i.e. $S_h \subseteq S_0$.
$\bullet$

\begin{Fo} \label{Sa6}
An identity between bracket symbols is valid
if and only if the identity is valid for $h=0$.
\end{Fo}

\begin{flushleft}
\bf  \large 4. Binary forms
\end{flushleft}

We use the algebra $H_I$ in order
to introduce noncommutative binary homogeneous $n$-forms.

For the following we suppose for the algebra $H_I$ that
$0\in I$.
We will use the short notations
\[  x:=x_0,\ \ \ \ \  y:=y_0.   \]
Furthermore let $H_{I \backslash \{ 0 \} }$
be the subalgebra of $H_I$
which is generated by the elements
$x_i,y_i$,   $i\in I \backslash \{ 0 \}$.

\begin{De} \label{De2} {\rm
We say that the nonvanishing element
\begin{eqnarray} \label{eq7} f:=
                                      x^n        A_0 +
       \left( \begin{array}{c} n \\ 1 \end{array}   \right) x^{n-1}y   A_1 +
       \left( \begin{array}{c} n \\ 2 \end{array}   \right) x^{n-2}y^2 A_2 +
 ....  + y^n A_n
\end{eqnarray}
of $H_I$ with $A_i\in H_{I\backslash \{0\} }$
is a {\it n-form },
if $f\in H_I^{Inv} $. }
\end{De}

\begin{Rem} \label{Rem4} \rm
The  representation of $f$  with coefficients
$A_i\in H_{I\backslash \{0\} }$ is unique.
This is a consequence of the PBW theorem.
\end{Rem}

\begin{Rem} \label{Rem5} \rm
Alternatively one can consider left coefficients.
For example, for
$f=(01)^2 =x^2A+2xyB+y^2C = A_Lx^2+2B_Lxy+C_Ly^2$
we have
\begin{eqnarray*}
A  =y_1^2,\ \ \ \ \  B  =-x_1y_1-\frac{h}{2}y_1^2,
\ \ \ \ \  C  =x_1^2+3hx_1y_1, \\
A_L=y_1^2,\ \ \ \ \  B_L=-x_1y_1+\frac{3h}{2}y_1^2,
\ \ \ \ \  C_L=x_1^2- hx_1y_1. \\
\end{eqnarray*}
\end{Rem}

\begin{Sa} \label{Sa7}
Let $f$ be an arbitrary $n$-form. Then the coefficients of $f$
obey the commutation relations of the coefficients
of the special $n$-form $(01)(02)...(0n)$.
\end{Sa}

\noindent {\it Proof:} Let
$f_1 =(0i_1)...(0i_n)=\sum {n\choose i} x^iy^{n-i}A_i$
and
$f_2=(0{j}_1)...(0{j}_n)=\sum {n\choose i} x^iy^{n-i}{A'}_i$
be two $n$-forms with $2n$ different indices $i_k,j_k$.
By (\ref{eq4}) we have commutation relations
of the form
\begin{eqnarray} \label{eq8}
{A'}_l {A}_k - {A}_k {A'}_l   =
\sum_{i\leq k,j\leq l,i+j<k+l}
       \alpha_{ij} {A}_i {A'}_j.
\end{eqnarray}
Because of the PBW property
the $\alpha_{ij}$ are uniquely determined.
These relations are also valid if not all $i_k,j_k$ are different.
In the case $f_1=f_2$ we can simplify the relations to
\begin{eqnarray} \label{eq9}
 {A}_l {A}_k - {A}_k {A}_l =
     \sum_{i\leq j,i\leq k,j\leq l,i+j<k+l}
     \beta_{ij} {A}_i {A}_j,
\end{eqnarray}
i.e. $\beta_{ij}=0$ for $i>j$.
According to Theorem \ref{Sa4}
two arbitrary $n$-forms $f$ and $f'$
have the representations
\[
\sum_{i_1,...,i_n\in I\backslash \{0\} }
(0i_1)(0i_2)...(0i_n) {P}_{i_1,...,i_n},\ \ \ \ \ \ \
\sum_{i_1,...,i_n\in I\backslash \{0\} }
(0i_1)(0i_2)...(0i_n) {P'}_{i_1,...,i_n}
\]
with central elements
${P}_{i_1,...,i_n},{P'}_{i_1,...,i_n}
\in H_{I\backslash \{0\} }^{Inv}$.
Therefore $f$ and $f'$ have the coefficients
\[
{A}_i =\sum A_i^{i_1,...,i_n} {P'}_{i_1,...,i_n},
\ \ \ \ \ \ \
{A'}_i=\sum A_i^{i_1,...,i_n} {P''}_{i_1,...,i_n}
\]
where $A_i^{i_1,...,i_n}$ are the coefficients of
the $n$-form $(0i_1)(0i_2)...(0i_n)$.
Therefore the equation (\ref{eq8}) is also valid for $f$ and $f'$.
It follows (\ref{eq9}) for $f=f'$.
$\bullet$

\begin{Rem} \label{Rem6} \rm
Proposition \ref{Sa7} is not true for the quantum group
$U_q (sl(2))$ (cf. \cite{Le}), i.e.
the commutation relations depend on the
concrete type of form.
\end{Rem}

\begin{Ex} \label{Ex3} \rm
Consider the quadratic form
\begin{eqnarray*}
f:=   (01)(02)=x^2A+xyB+y^2C
\end{eqnarray*}
with
\begin{eqnarray*}
A  =y_1y_2,
B  =-\frac{h}{2}(x_1y_2+y_1x_2)-hy_1y_2,
C  =x_1x_2+hx_1y_2+2hy_1x_2+2h^2y_1y_2. \\
\end{eqnarray*}
For $f$ we have the following commutation relations
\begin{eqnarray} \label{eqz1}
\begin{array}{lll}
CA & = &  AC - 4h\ AB + 6h^2\ A^2,           \\
CB & = &  BC - 2h\ AC + 2h^2\ AB -3h^3\ A^2, \\
BA & = &  AB - 2h\ A^2,
\end{array}
\end{eqnarray}
\begin{eqnarray*}
\begin{array}{llllll}
Ax & = &  xA-2hyA,
\ \ \ \ \ \ \ & Ay & = &  yA,  \\
Bx & = &  xB+hxA+h^2yA,
\ \ \ \ \ \ \ & By & = &  yB-hyA,  \\
Cx & = &  xC+2hyC+2hxB+4h^2yB,
\ \ \ \ \ \ \ & Cy & = &  yC-2hyB+2h^2yA.
\end{array}
\end{eqnarray*}
\end{Ex}

\begin{Ex} \label{Ex4} \rm
For the cubic form
\begin{eqnarray*}
f=(01)(02)(03)
\end{eqnarray*}
we have
\[  \begin{array}{l}
A=  y_1y_2y_3, \\
B=
 -\frac{h}{3}(x_1y_2y_3+y_1x_2y_3+y_1y_2x_3)-hy_1y_2y_3, \\
C=\frac{h}{3}(x_1x_2y_3+x_1y_2x_3+y_1x_2x_3)
 +\frac{2h}{3}x_1y_2y_3
 +hy_1x_2y_3
 +\frac{4h}{3}y_1y_2x_3
 +2h^2y_1y_2y_3, \\
D=  -x_1x_2x_3-hx_1x_2y_3-2hx_1y_2x_3-3hy_1x_2x_3-2h^2x_1y_2y_3
-3h^2y_1x_2y_3-6h^2y_1y_2x_3-6h_3y_1y_2y_3
\end{array}  \]
and the commutation relations
\begin{eqnarray*}
\begin{array}{ll}
 DA= &  AD- 9 h         AC
                               +36 h^2 AB
                               -60 h^3 A^2, \\
 DB= & BD      -3 h AD
                               -3 h BC
                             +9 h^2 AC
                             +6 h^2 B^2
                            -24 h^3 AB
                            +36 h^4 A^2, \\
 DC= & CD          -6h BD
                               +3h C^2
                           +12 h^2 AD
                            -6 h^3 AC
                           -18 h^3 B^2
                           +36 h^4 AB
                           -48 h^5 A^2, \\
 CA= & AC            -6h AB
                           +12 h^2 A^2, \\
 CB= & BC            - h AC
                              -2 h B^2
                            +2 h^2 AB
                            -4 h^3 A^2, \\
BA= &
AB            -3h A^2.
\end{array}
\end{eqnarray*}
\end{Ex}

\begin{flushleft}
\bf \large 5. Linear forms and points
\end{flushleft}

\noindent
The simplest example of a linear form is
\begin{eqnarray*}
xA+yB:=(01)=xy_1-y(x_1+hy_1),
\end{eqnarray*}
i.e. $A=y_1$, $B=-x_1-hy_1$.

Another example is the form
\begin{eqnarray*}
xA+yB:=(02)(13)+(03)(12)
\end{eqnarray*}
corresponding to the fourth harmonic point
with
\begin{eqnarray*}
\begin{array}{l}
A =  2 x_1y_2y_3- y_1x_2y_3- y_1y_2x_3-3hy_1y_2y_3 \\
B = -  x_1x_2y_3- x_1y_2x_3+2y_1x_2x_3- hx_1y_2y_3
  +hy_1x_2y_3+3hy_1y_2x_3+3h^2y_1y_2y_3.
\end{array}
\end{eqnarray*}
\noindent
Let $f=xA+yB$ be a linear form and let
\begin{eqnarray*}
Y=A, \ \ \ \ \ \ \ X=-B-hA.
\end{eqnarray*}
Then $f$ has the representation
\begin{eqnarray*}
f=xY-yX-hyY=-(Xy-Yx-hYy).
\end{eqnarray*}

\begin{Sa} \label{Sa8}
Let $f=xY-yX-hyY$ and
$\overline{f}=xY'-yX'-hyY'$ be two linear forms.
The elements $X,Y,X',Y'$ have
commutation relations similar to (\ref{eq4}),
i.e. we have
\begin{eqnarray} \label{eq10}
\begin{array}{ll}
XY  = YX+hY^2,\ \ \ \ \ \ \ &
X'Y'= Y'X'+h{Y'}^2, \\
XX' = X'X-hX'Y+hY'X+h^2Y'Y,\ \ \ \ \ \ \   &
YY' = Y'Y, \\
XY' = Y'X+hY'Y,\ \ \ \ \ \ \ &
YX' = X'Y-hY'Y.
\end{array}
\end{eqnarray}
\end{Sa}

\noindent
The proof is analogous to Proposition \ref{Sa7}.

\begin{Rem} \label{Rem7} \rm
We consider $X,Y$ as homogeneous coordinates
of a noncommutative point
with the projective coordinate $Z:=X\ Y^{-1}+\frac{h}{2} \in
Q_I$.
\end{Rem}

\begin{flushleft}
\bf \large 6. Invariants of Forms
\end{flushleft}

We consider a
$n$-form $f=\sum_{i=0}^n
{n\choose i}
x^{n-i} y^i A_i $.

\begin{De} \label{De3}  {\rm
We say that the polynomial expression
of the coefficients of $f$
\[ I_f = I_f (A_0,...,A_n) \]
is an {\it invariant},
if $I_f \in H_{I \backslash \{ 0 \} }^{Inv}$.
We say that the polynomial expression
\[ C_f = C_f ( A_0,...,A_n ,x_0,y_0 ) \]
is a  {\it covariant},
if $C_f$ is a $m$-form.}
\end{De}

\begin{Ex} \label{Ex5} \rm
The simplest invariant of a quadratic form is the discriminant \\
$AC - B^2 - h AB + \frac{3h^2}{4} A^2$ (cf. below).
\end{Ex}

\begin{flushleft}
\bf \large 7. The Symbolic Method
\end{flushleft}

We consider the $n$-form
\[ { f} =                                          x^n        { A}_0 +
\left( \begin{array}{c} n \\ 1 \end{array} \right) x^{n-1}y   { A}_1 +
\left( \begin{array}{c} n \\ 2 \end{array} \right) x^{n-2}y^2 { A}_2 +
                                      ....       + y^n        { A}_n. \]
Our aim is to construct invariants of degree $k$.

We consider the collection of $k$ very simple $n$-forms
\[ f':=(01)^n,\ \ f'':=(02)^n,\ \  ... ,\ \  f^{(k)}:=(0k)^n. \]
Because of (\ref{eq4}) the form $f^{(j)}$ has the
coefficients $A^{(j)}_i$ with
\begin{eqnarray*}
A^{(j)}_i =(-1)^{i} { n \choose i } x_j^{i}y_j^{n-i}+
            \sum_{k=1}^{i} c_{ik} h^{k} x_j^{i-k}y_j^{n-i+k}
\end{eqnarray*}
and certain numbers $c_{ik}\in \ZZ$.

Consider an element $d$ of $H^{Inv}_I$
which is homogeneous of degree $n$
with respect to the indices $1,2,...,k$.
(For example $(12)(13)(23)$ is
homogeneous of degree 2 with respect to $1, 2, 3$.)

Because of the PBW theorem for $H_I$
and because the application of the reduction rules (\ref{eq4})
does not change the degree of
homogeneity with respect to $1,2,...,k$
we can express $d$ as a linear combination of monomials
\begin{eqnarray} \label{eq11}
   x_1^{ n- i_1 } y_1^{ i_1 }
   x_2^{ n- i_2 } y_2^{ i_2 } ...
   x_k^{ n- i_k } y_k^{ i_k }.
\end{eqnarray}

In these products we replace expressions $x_j^n$ by
\begin{eqnarray*}
=(-1)^{n} (A^{(j)}_n - \sum_{k=1}^{n} c_{nk} h^{k}
x_j^{n-k}y_j^{k}).
\end{eqnarray*}
We obtain monomials of products of terms
$x_j^{n-i}y_j^{i}$ with $1\leq i\leq n$.
Then we replace the products $x_j^{n-1}y_j^{1}$ by
\begin{eqnarray*}
=(-1)^{n-1} { n \choose 1 }^{-1}
(A^{(j)}_{n-1} - \sum_{k=1}^{n-1} c_{n-1,k} h^{k}
x_j^{n-1-k}y_j^{1+k}),
\end{eqnarray*}
etc..
In the $n+1$-th step we replace
the terms $y_j^n$ by $A^{(j)}_0$.

In this way
we obtain for the symbol $d$ the representation
\begin{eqnarray*}
d = \sum c_{i_1,i_2,...,i_k}
A^{(1)}_{i_1} A^{(2)}_{i_2} ... A^{(k)}_{i_k}.
\end{eqnarray*}

Furthermore we consider the elements
\[    {\cal I}' = \sum c_{i_1,i_2,...,i_k}
{\cal A}^{(1)}_{i_1} {\cal A}^{(2)}_{i_2} ... {\cal A}^{(k)}_{i_k} \]
for $k$ arbitrary $n$-forms $F^{j}$
with coefficients ${\cal A}^{(j)}_i$.

We form
\[ E({\cal I}') =\sum  d_{i_1,i_2,...,i_k}
{\cal A}^{(1)}_{i_1} {\cal A}^{(2)}_{i_2} ... {\cal A}^{(k)}_{i_k}
\]
and consider the homomorphism $h_{f',...,f^{(k)} }$
determined by
\[ h_{f',...,f^{(k)}} ({\cal A}_i^{j} ):= A_i^{j}.  \]
It follows
\[   0 = E(d) = E(h_{f',...,f^{(k)} } ({\cal I}')) =
                h_{f',...,f^{(k)} } (E({\cal I}')) =
     \sum d_{i_1,i_2,...,i_k}
     A^{(1)}_{i_1} A^{(2)}_{i_2} ... A^{(k)}_{i_k}    \]
\[ =\sum_{i_1,...,i_k} d_{i_1,i_2,...,i_k} (-1)^{i_1+...+i_k}
\left( {n \choose i_1}  ... {n \choose i_k} \right)^{-1}
            x_1^{ i_1 } y_1^{ n - i_1 }
            x_2^{ i_2 } y_2^{ n - i_2 } ...
            x_k^{ i_k } y_k^{ n - i_k }     \]

\[ +\sum_{i_1,...,i_k} d_{i_1,i_2,...,i_k} \sum_{1\leq l_1 \leq i_1,
                                                     ...,
                                                 1\leq l_k \leq i_k      }
              c_{i_1 l_1}...c_{i_k l_k} h^{l_1+...+l_k}
            x_1^{ i_1-l_1 } y_1^{ n-i_1+l_1 }
            x_2^{ i_2-l_2 } y_2^{ n-i_2+l_2 } ...
            x_k^{ i_k-l_k } y_k^{ n-i_k+l_k }     \]
The PBW theorem allows a comparison of coefficients.
For the coefficient at
           $x_1^{ i_1 } y_1^{ n - i_1 }
            x_2^{ i_2 } y_2^{ n - i_2 } ...
            x_k^{ i_k } y_k^{ n - i_k } $ we have
\[ 0 =  d_{i_1,i_2,...,i_k} (-1)^{i_1+...+i_k}
\left( {n \choose i_1}  ... {n \choose i_k} \right)^{-1}
+\]
\[       +  \sum_{1\leq l_1 \leq n-i_1,..., 1\leq l_k \leq n-i_k  }
              d_{i_1+l_1,...,i_k+l_k}
              c_{i_1+l_1, l_1}...c_{i_k+l_k l_k} h^{l_1+...+l_k}
\]
These equations form a triangular equation system
for the elements $d_{i_1,...,i_k}$
with the unique solution $d_{i_1,...,i_k}=0$,
i.e. $E( {\cal I}' )=0$.
Similarly we obtain $F( {\cal I}' )=0$ and
$H( {\cal I}' )= 0$.
Therefore $ {\cal I}' $ is a common
invariant of $k$ arbitrary $n$-forms
${ F}'$, ${ F}''$, ... , ${ F}^{(k)}$.

In particular we obtain an invariant of $f$,
if we replace ${ F}'$, ${ F}''$, ... , ${ F}^{(k)}$
by $f$. Therefore
\[ {\cal I} = \sum c_{i_1,i_2,...,i_k}
   { A}_{i_1} { A}_{i_2} ... { A}_{i_k}  \]
is an invariant of $f$.

Therefore we have demonstrated the following Theorem.

\begin{Th} \label{Sa9}
Let $f$ be a $n$-form
and let $d$ be an element of $H_I^{Inv}$
which is homogeneous of degree $n$ with respect to $1,2,...,k$.
Then the above construction gives
an invariant $\cal I$ of $f$ of
degree $k$.
\end{Th}

Analogously, the symbolic method admits the
construction of covariants, common invariants
and common covariants.

The simplest symbols are
$(12)^2$ and  $(12)(13)(23)$
for the quadratic form,
$(12)^3$
for the cubic form and
$(12)^4$, $(12)^2 (13)^2 (23)^2$
for the quartic form.

\begin{Ex} \label{Ex6} \rm
We apply the symbolic method in the simplest example.\\ Let $f =
x^2A + 2xyB + y^2C$ be a quadratic form. For the symbol $(12)^2$
we obtain
\[ (12)^2=
   x_1^2  y_2^2
 - 2x_1 y_1 x_2 y_2
 + y_1^2  x_2^2
 - hx_1 y_1 y_2^2
 + 3hy_1^2  x_2 y_2
  \]
\[ =
   (C^{(1)}-3hx_1y_1) y_2^2
 - 2x_1 y_1 x_2 y_2
 + y_1^2  (C^{(2)}-3hx_2y_2)
 - hx_1 y_1 y_2^2
 + 3hy_1^2  x_2 y_2  \]

\[ =
   \left(C^{(1)}-3h(-B^{(1)}-\frac{h}{2}y_1^2)\right) y_2^2
\ -\  2(-B^{(1)}-\frac{h}{2}y_1^2) (-B^{(2)}-\frac{h}{2}y_2^2)\]
\[ + y_1^2  \left(C^{(2)}-3h(-B^{(2)}-\frac{h}{2}y_2^2)\right)
  \ -\ h(-B^{(1)}-\frac{h}{2}y_1^2) y_2^2
 \ +\ 3hy_1^2  (-B^{(2)}-\frac{h}{2}y_2^2)  \]

\[ =
   \left(C^{(1)}-3h(-B^{(1)}-\frac{h}{2}A^{(1)})\right) A^{(2)}
\ -\  2(-B^{(1)}-\frac{h}{2}A^{(1)}) (-B^{(2)}-\frac{h}{2}A^{(2)})\]
\[ + A^{(1)}  \left(C^{(2)}-3h(-B^{(2)}-\frac{h}{2}A^{(2)})\right)
  \ -\ h(-B^{(1)}-\frac{h}{2}A^{(1)}) A^{(2)}
 \ +\ 3hA^{(1)}  (-B^{(2)}-\frac{h}{2}A^{(2)})  \]
\begin{eqnarray*}
=A^{(1)}C^{(2)}+C^{(1)}A^{(2)}-2   B^{(1)}B^{(2)}
       +3h  B^{(1)}A^{(2)}
       -h   A^{(1)}B^{(2)}
       +\frac{3h^2}{2}A^{(1)}A^{(2)}.
\end{eqnarray*}
We identify the coefficients of $(01)^2$ and $(02)^2$
and obtain the discriminant
\begin{eqnarray*}
d_2=AC+CA-2   B^2
       +3h  BA
       -h   AB
       +\frac{3h^2}{2}A^2.
\end{eqnarray*}
We can apply the commutation relations (\ref{eqz1})
in order to get the ordered expression
\begin{eqnarray*}
d_2=2 AC - 2 B^2 - 2h AB + \frac{3h^2}{2} A^2
\end{eqnarray*}
or the symmetric expression
\begin{eqnarray*}
d_2=AC+CA-2 B^2+h(AB+BA)-\frac{5h^2}{2}A^2
\end{eqnarray*}
or expressions without the term $AB$ or $A^2$
\begin{eqnarray*}
d_2=\frac{3}{2} AC+\frac{1}{2}CA-2B^2-\frac{3h^2}{2}A^2
\end{eqnarray*}
\begin{eqnarray*}
= AC+ CA -2 B^2- \frac{h}{4}  AB + \frac{9h}{4} BA.
\end{eqnarray*}
For the special forms $f=(01)(02)$ and $f=(01)^2$
we obtain
$d_2=-\frac{1}{2}(12)^2$ and $d_2=0$, respectively.
$(12)(13)(23)$ gives the zero invariant.
\end{Ex}

\begin{Ex} \label{Ex7} \rm
Consider the cubic form $f$.
The simplest symbol $(12)^3$ gives the zero invariant.
The Hessian covariant  $\Delta$
has the symbol $(01)(02)(12)^2$.
Up to a constant factor we obtain
\[ \Delta = x^2 K + xy L + y^2 M \]
with
\begin{eqnarray*}
\begin{array}{ll}
K= &     2 AC
          -2B^2
        -4 hAB
     +2 h^2 A^2, \\
L= &   2   AD
         -2 BC
        -2h AC
     +4 h^2 AB
     -2 h^3 A^2,  \\
M=  &  2    BD
          -2 C^2
        -4 h AD
        +6h^2B^2
       -6h^3 AB
       +8h^4 A^2.
\end{array}
\end{eqnarray*}
The discriminant $d_3$ of $f$
has the symbol $(12)^2(34)^2(13)(24)$.
Up to a scalar factor $d_3$ is also the discriminant
of the quadratic form $\Delta$.
We obtain
\begin{eqnarray*}
2 d_3=2 KM - \frac{1}{2} L^2 - {h} KL + \frac{3h^2}{2} K^2
\end{eqnarray*}
and
\begin{eqnarray*}
\begin{array}{lc} d_3 =
 &          - A^2 D^2
         + 6\ A B C D
         - 4\ A C^3
         - 4\ B^3 D
         + 3\ B^2 C^2 \\
& +9h    (-2\ A B^2 D
           +\ A B C^2
           -\ A^2 C D)
+h^2( -    9\ B^4
      +   40\ A^2 B D
         - 7\ A^2 C^2
        + 12\ A B^2 C) \\
& +6h^3(- 12\ A^3 D
      - 6\    A^2 B C
    -   11\   A B^3)
 +2h^4( 75\   A^2 B^2
     +  38\   A^3 C)\\
&  - 384\ h^5 A^3 B
   + 652\ h^6 A^4.
\end{array}
\end{eqnarray*}
Furthermore the cubic covariant $j$
has the symbol $(01)^2(02)(13)(23)^2$.
We have
\begin{eqnarray*}
\begin{array}{ll}
j= &
\ \  x^3(
           A^2 D -
         3 A B C +
         2 B^3   +
        9h A B^2 -
      2h^2 A^2 B +
      6h^3 A^3   )   \\
& +x^2y(+
          3  A B D -
          6  A C^2 +
          3  B^2 C-
         9h  A^2 D+
        12h  A B C -
         3h  B^3-
       6h^2  A B^2-
      12h^3  A^2 B
     + 6h^4  A^3    )   \\
& +xy^2(-
           3  A C D+
           6  B^2 D-
           3  B C^2 -
          6h  A B D+
          9h  A C^2-
          3h  B^2 C
\\ & \ \ \ \ \ \
          +
       30h^2  A^2 D-
       12h^2  A B C+
       12h^2  B^3-
       12h^3  A B^2+
       60h^4  A^2 B
      -90h^5  A^3
         )  \\
& +y^3(- A D^2+
     3  B C D-
     2  C^3-
   15h  B^2 D+
    6h  B C^2+
 36h^2  A B D -
  8h^2  A C^2+
  6h^2  B^2 C
\\ & \ \ \ \ \ \
-102h^3  A^2 D-
 18h^3  A B C  -
 42h^3  B^3  +
 38h^4  A^2 C  +
126h^4  A B^2 -
378h^5  A^2 B
+704h^6  A^3).
\end{array}
\end{eqnarray*}
\end{Ex}

\begin{Rem} \label{Rem8} \rm
In the classical case, the zeros of $j$ are the three fourth
harmonic points of the zeros of $f$ and the zeros of $\Delta$ are
the two equiharmonic points of the zeros of $f$ or of $j$
(cf. \cite{Cl}).
\end{Rem}

\begin{Rem} \label{Rem9} \rm
For $f,\Delta,j,d_3$ we have the syzygy
\begin{eqnarray} \label{eq12}
4\ \Delta^3\  +\  j^2\ +\ d_3\ f^2\ =\ 0.
\end{eqnarray}
By Theorem 4
we can represent $f,\Delta,j,d_3$ by bracket symbols.
The coefficients do not depend on $h$.
By Corollary \ref{Sa6} the equation follows
from the classical situation, cf. \cite{Cl,Ol}.
\end{Rem}

\begin{flushleft}
\bf \large 8. Polarization
\end{flushleft}

In this section we generalize the classical polarization operators
$x_l \frac{\partial}{\partial x_k}+
 y_l \frac{\partial}{\partial y_k}$ by
operators $\Delta_{kl}: H_I \rightarrow  H_I $ with
$k,l\in I$.
We fix a PBW basis, which consists of the elements
$    x_{a_1}^{ i_1 }   y_{a_1}^{ j_1 }
     x_{a_2}^{ i_2 }   y_{a_2}^{ j_2 }  ...
     x_{a_k}^{ i_k }   y_{a_k}^{ j_k } $
with $a_1 < a_2 < ... < a_k$. We set
\begin{eqnarray*}
\Delta_{kl} ( x_k ) = x_l,\ \ \ \ \ \ \ & \Delta_{kl} ( y_k ) =
y_l
\end{eqnarray*}
and
\begin{eqnarray*}
\Delta_{kl} ( x_i ) = 0,  \ \ \ \ \ & \Delta_{kl} ( y_i ) = 0
\end{eqnarray*}
for $i\neq k$ and require that
\[ \Delta_{kl} ( ab ) =     \Delta_{kl} (a)\ b +
                         a\ \Delta_{kl} (b)        \]
for $a,b\in H_I$.

\begin{Sa} \label{Sa10}
$\Delta_{kl}$ is a well defined linear operator.
\end{Sa}

\noindent
{\it Proof:}
One has to check the consistency of the Leibniz rule
with the relations (\ref{eq4}).
For example we have
\[\Delta_{kl}(y_ix_k) = y_i x_l = x_l y_i - h y_l y_i =
  \Delta_{kl}(x_k y_i - h y_k y_i).\bullet  \]

Furthermore we define the operator
$P_{ij}: H_I\rightarrow H_I$ with
\[ P_{ij}a={1 \over n}\Delta_{ij}(a) \]
if $a$ is homogeneous of degree $n$
with respect to the index $i$.
Analogous we define polarization operators
$\Delta_{(x_i,y_i),(X,Y)}$ and
$P_{(x_i,y_i),(X,Y)}$
where $X,Y$ are the homogeneous coordinates
of a linear form.
$P_{ij}$ and $P_{(x_i,y_i),(X,Y)}$
have the simple property
\begin{eqnarray} \label{eq13}
\begin{array}{llllll}
  P_{(x_i,y_i),(X,Y)} (ik) & = & (Xk),\ \ \ \ \ \ \ &
  P_{(x_i,y_i),(X,Y)} (kl) & = &  0
\end{array}
\end{eqnarray}
for  $i\neq k,l$
with $(Xk) := Xy_k-Yx_k-hYy_k \in H_I^{Inv}$.
In particular we have
\begin{eqnarray} \label{eq13d}
\begin{array}{llllll}
  P_{ij} (ik) & = & (jk),  \ \ \ \ \ \ \ &
  P_{ij} (kl) & = &  0
\end{array}
\end{eqnarray}
for  $i\neq k,l$.
We have the Leibniz rule
\begin{eqnarray} \label{eq14}
\begin{array}{rcl}
P_{(x_i,y_i),(X,Y)} ( ab ) & = &
n P_{(x_i,y_i),(X,Y)} (a)\ b +m a\ P_{(x_i,y_i),(X,Y)} (b)
\end{array}
\end{eqnarray}
for $a,b \in H_I$ with $\deg_i(a)=m$ and $\deg_i(b)=n$.
It follows
\[   P_{(x_i,y_i),(X,Y)} H_I^{Inv} \subseteq H_I^{Inv}.     \]

\begin{Rem} \rm \label{Rema1}
Let $f$ be an element of $H_I$,
let $xY-yX-hyY$ be a linear form
and let $G=G_{(x_i,y_i),(X,Y)}$
be the algebra homomorphism which is determined by
$G(x_i)=X$, $G(y_i)=Y$ and
$G(x_j)=x_j$, $G(y_j)=y_j$ for $j\neq i$.
Because of Proposition \ref{Sa8} $G$ is compatible with (5).
Therefore the algebra homomorphism $G$ is well defined.
We can extend $G$ to $Q_I$.
\end{Rem}

\begin{Le} \label{Sa11}
Let $f$ be a $k$-form
and let $xY-yX-hyY$ be a linear form.
If we replace $x,y$ in $f$ by $X,Y$, respectively,
we obtain $P_{(x,y),(X,Y)}^k f$,
i.e. we have $P_{(x,y),(X,Y)}^k f = G_{(x,y),(X,Y)} f$.
\end{Le}
\noindent {\it Proof:} The Lemma is true in the classical case
(cf. \cite{Cl}). We represent $f$ by brackets $(ij)$. By
(\ref{eq13}), (\ref{eq14}) the action of $P_{(x_i,y_i),(X,Y)}$ on
the algebra of brackets is equal to the classical case. Therefore
we obtain $P_{(x_i,y_i),(X,Y)}^k f$ if we replace all $(0i)$ by
$(Xi)$. By Remark \ref{Rema1} the proposition follows.$\bullet$

\begin{Rem} \label{Rem10} \rm
An analogue of Lemma \ref{Sa11} is not true in the case of
$U_q(sl(2))$-symmetry. For the noncommutative coordinate algebra
in \cite{Le} we can define polarisation operators but we have no
analogue of $G$.
\end{Rem}

\begin{Le} \label{Sa12}
A formula which contains only brackets $(ij)$, polarization
operators $P_{(x_k,y_k),(X,Y)}$
and complex coefficients is true if and only if
this formula is valid in the classical case $h=0$.
\end{Le}

\noindent
{\it Proof:}
The proof follows from (\ref{eq13}), (\ref{eq14})
and Corollary \ref{Sa6}.$\bullet$

\begin{Ex} \label{Ex8} \rm
By the corresponding classical result
(cf. \cite{Cl})
and Lemma 12 we obtain that
$P_{01}^{n-1} f$ is a linear form with
$P_{01}^{n-1} f = P_{10} P_{01}^{n} f$.
\end{Ex}

\begin{flushleft}
\bf \large 9. Typical representation
\end{flushleft}

We will use the typical representation of a binary form
(cf. \cite{Cl,He,Ko} for the classical case)
for the decomposition into linear forms.
We consider the n-form
\begin{eqnarray} \label{eq15}
f = x^n A_0 + x^{n-1} y A_1 + x^{n-2} y^2 A_2+
    \cdots +y^n A_n
\end{eqnarray}
with $A_i \in H_{I\backslash \{0,1\}}$
and the $n$-th polar
\[  f_1 :=  P_{01}^n f =
 x_1^n A_0 + x_1^{n-1} y_1 A_1 + x_1^{n-2} y_1^2 A_2+
    \cdots +y_1^nA_n,              \]
cf. Lemma \ref{Sa11}.
Furthermore let $\xi,\eta$ be the two linear forms
\begin{eqnarray*}
\begin{array}{cccccc}
\xi  & = &     P_{10} f_1 &= P_{01}^{n-1} f  & =: &
x \delta - y \gamma -h y \delta
\\
\eta & = &     (01)  &      & = &  xy_1-yx_1-hyy_1
\end{array}
\end{eqnarray*}
with $\gamma,\delta \in H_{I\backslash \{0\} }$.
$\gamma,\delta$
have commutation relations analogous to $x_i,y_i$.

\begin{Th} \label{Sa13}
We have
\begin{eqnarray} \label{eq16}
 f\ =  f_1^{-n}
 ( \xi^n u_0 + \xi^{n-1} \eta u_1 + \xi^{n-2} \eta^2 u_2+
           \cdots +\eta^nu_n )
\end{eqnarray}
with $u_i\in H_{I\backslash \{0\}}$.
The coefficients $u_0,u_1,\cdots,u_n$ are covariants
with respect to $x_1,y_1$.
\end{Th}

\noindent {\it Proof}:
We refer to \cite{Cl} for the classical case.
The Theorem follows from Lemma \ref{Sa12}.$\bullet$

\begin{Rem} \label{Rem11} \rm
We call (\ref{eq16}) the {\it typical represetation} of $f$.
We mention that $u_0=f_1$ and $u_1=0$ (cf. \cite{Cl}).
One can derive the typical representation
if one expresses $x,y$ by $\xi,\eta$. We have
\begin{eqnarray*}
\begin{array}{ccc}
   \xi  y_1 &  = &  x  \delta  y_1 - y (\gamma + h \delta ) y_1,
   \\
   \eta  \delta  &  = &  x y_1  \delta  - y (x_1+hy_1)  \delta.
\end{array}
\end{eqnarray*}
Because of $ y_1 \delta=\delta y_1 $ it follows
\begin{eqnarray} \label{eq17}
  \xi y_1 -  \eta \delta = y ( x_1 \delta - \gamma y_1 ).
\end{eqnarray}
Furthermore we have
\begin{eqnarray*}
\begin{array}{ccc}
   \xi  x_1 &  = &  x  \delta  x_1 - y (\gamma + h \delta ) x_1,
   \\
   \eta  \gamma  &  = &  x y_1  \gamma  - y (x_1+hy_1)  \gamma .
\end{array}
\end{eqnarray*}
Because of $(x_1+h y_1) \gamma = (\gamma+h \delta) x_1$
we have
\begin{eqnarray} \label{eq18}
   \xi x_1 - \eta \gamma  =    x ( \delta x_1 -  y_1 \gamma ).
\end{eqnarray}
We mention the relation $ \delta x_1 - y_1 \gamma   =
        x_1 \delta - \gamma y_1   =
        x_1 \delta -y_1 \gamma -h y_1 \delta =  f_1$.
From (\ref{eq17}) and (\ref{eq18}) we obtain the inverse transformation
\begin{eqnarray} \label{eq19}
\begin{array}{l}
   x   =  ( \xi x_1 - \eta \gamma ) f_1^{-1}, \\
   y   =  ( \xi y_1 - \eta \delta ) f_1^{-1}.
\end{array}
\end{eqnarray}
We insert the expressions for $x,y$ into (\ref{eq15})
and obtain
\begin{eqnarray*}
\begin{array}{c}
 f = f_1^{-n}
   (  (\xi x_1 - \eta \gamma  )^n
                         A_0 +
   (   \xi x_1 - \eta \gamma  )^{n-1}
    ( \xi y_1 - \eta \delta ) A_1 +
           \cdots +
    ( \xi y_1 - \eta \delta )^n
          A_n).
\end{array}
\end{eqnarray*}
The elements $\xi,\eta$ are central.
It follows
\[ f\ =  f_1^{-n}
 ( \xi^n c_0 + \xi^{n-1} \eta c_1 + \xi^{n-2} \eta^2 c_2+
           \cdots +\eta^n c_n )  \]
with $c_i\in H_{I\backslash \{0\}}$.
It is possible to show that $c_i=u_i$.
\end{Rem}

\begin{flushleft}
\bf \large 10. The decomposition of binary forms into linear forms
\end{flushleft}

By (\ref{eq16})
we can consider $f$ as a commutative $n$-form with
variables $\xi,\eta \in Q_I^{Inv}$ and
coefficients $u_i \in Q_{ I\backslash \{0\} }^{Inv} $.
It follows from Corollary \ref{Sa6}
that there is no algebraic relation between
$\xi$, $\eta$ and $u_i$, $i=1,...,n$.
Therefore there is a commutative splitting field
$\Sigma \supset Q_{I\backslash \{0\} }^{Inv}$
with
\begin{eqnarray} \label{eq20}
 f = C\ (\xi-\sigma_1\eta)(\xi-\sigma_2\eta)
            \cdots (\xi-\sigma_n\eta)
\end{eqnarray}
where $\sigma_i\in \Sigma $ and
$C = u_0 f_1^{-n} = f_1^{1-n}
\in Q^{Inv}_{I\backslash \{0\}}$.

We can suppose, that the elements
$\sigma_1,...,\sigma_n$
are invariant
with respect to an extended action of $U_h(sl(2))$
and that they commute with the skew field $Q_I$.

\begin{Rem} \label{Rem12} \rm
In the classical case the elements $\sigma_i$ are the associated
roots of G. Kohn (cf. \cite{Ko}).
\end{Rem}

We consider the central skew field extensions
$Q_{I\backslash \{0\} }^\Sigma:=
 Q_{I\backslash \{0\} }[\sigma_1,...,\sigma_n]$
and \\
$Q^\Sigma_I:=Q_I[\sigma_1,...,\sigma_n]$, (cf. \cite{Co}).
The elements $(\xi-\sigma_i\eta )\in Q^\Sigma_I$
are linear forms in $x,y$. We have
\begin{eqnarray*}
\begin{array}{lll}
\xi - \sigma_i \eta & = & x(\delta -\sigma_i y_1 )-
                          y(\gamma -\sigma_i x_1)-
                         hy(\delta -\sigma_i y_1 )
=: x Y_i-y X_i-hy Y_i
\end{array}
\end{eqnarray*}
with
$X_i=\gamma -\sigma_i x_1 $ and
$Y_i=\delta -\sigma_i y_1 \in Q_{I\backslash \{0\} }^\Sigma$.
Therefore we have
\begin{Th} \label{Sa14} Let
$f=\sum_{i=1}^n x^{n-i}y^iA_i$ be a binary form
with  $A_i\in H_{I\backslash \{0,1\} }$.
Then we have the decomposition
\begin{eqnarray} \label{eq21}
f= C\ \prod_{i=1}^n (xY_i-yX_i-hyY_i)
\end{eqnarray}
into (commuting, invariant) linear forms
where $X_i,Y_i$ are elements of the
central skew field extension $Q_{I\backslash \{0\} }^\Sigma  $
and $C= f_1^{1-n}
\in Q^{Inv}_{I\backslash \{0\} }$.
\end{Th}

We can consider $X_i,Y_i$ as homogeneous coordinates
of the unique $n$ zeros \\
$Z_i=X_iY_i^{-1}+\frac{h}{2}$
of $f$ in the following sense.
\begin{Sa} \label{Sa15}
Let $f$ be a $n$-form with the decomposition (\ref{eq21}).\\
 (i) We have $P_{(x,y),(X_i,Y_i)}^n f = 0$.\\
(ii) Let $xY-yX-hyY$ be a linear form with
     $P_{(x,y),(X,Y)}^n f = 0$.
     Then there is an $i\in \{1,2,...,n\}$ with
     $ Z_i = X Y^{-1} + \frac{h}{2} $.
\end{Sa}

\noindent
{\it Proof:}
(i) follows from Lemma \ref{Sa11}
and $X_iY_i-Y_iX_i-hY_i^2=0$.  \\
(ii) Let
\[    P_{(x,y),(X,Y)}^n f =
C\ \prod_{i=1}^n (XY_i-YX_i-hYY_i)=0. \]
Because $Q_{I}^\Sigma  $
has no zero divisors there is an $i$ with $XY_i-YX_i-hYY_i=0$.
It follows
$( X - h Y ) Y_i = Y X_i$    $\Longrightarrow $
$ Y^{-1}X -h = X_iY_i^{-1} $ $\Longrightarrow $
$ XY^{-1} = X_iY_i^{-1}$ $\Longrightarrow $
$ XY^{-1}+ \frac{h}{2} = X_iY_i^{-1}+ \frac{h}{2}$.
$\bullet$

\begin{Rem} \label{Rem13} \rm
$X_i,Y_i$ and $C$ contain the additional
variables $x_1,y_1$ which are
not contained in the coefficients $A_i$, $i=0,...,n$. It is
possible to replace $x_1,y_1$ by the coefficients $X,Y$ of an
arbitrary linear form with $P_{(x,y),(X,Y)} f \neq 0$. Then the
analogous construction gives the points ${X'}_i,{Y'}_i$ with
${X'}_i {Y'}_i^{-1} ={X}_i {Y}_i^{-1}  $.
Therefore $Z_i$ is independent of $x_1,y_1$.
\end{Rem}

In the following sections we will decompose
the forms with $n=2,3$.

\begin{flushleft}
\bf \large 11. The quadratic equation
\end{flushleft}

At first we consider the classical case.
The classical quadratic form $f=ax^2+2bxy+cy^2$ has
the decomposition
\[   f = \frac{1}{ P_{01}^2 f } (xy_+ -yx_+) (xy_- -yx_-)   \]
with
\[   x_\pm =  -bx_1-cy_1\pm x_1\sqrt{b^2-ac}
\ \ \ \ \ \ \ \ \ \ {\rm and}
\ \ \ \ \ \ \ \ \ \   y_\pm =   ax_1+by_1\pm y_1\sqrt{b^2-ac}.    \]
This decomposition corresponds to the solution
of the quadratic equation $az^2+2bz+c=0$ of Clebsch
\begin{eqnarray*}
 z_\pm =  \frac{-bz_1-c \pm z_1\sqrt{b^2-ac} }{
                   az_1+b \pm    \sqrt{b^2-ac} }
\end{eqnarray*}
(cf. \cite{Cl,Mt}).

Now let
\[  f = x^2A+2xyB+y^2C           \]
be a noncommutative quadratic form
with $A,B,C \in H_{I\backslash \{0,1\}}$
with the discriminant
\[
   d_2=2 AC - 2 B^2 - 2h AB + \frac{3h^2}{2} A^2.
\]
We have the equation
\[  f \ P_{01}^2 f\ =\ (P_{01}f)^2\ + \frac{1}{2}\ d_2\ (01)^2 \]
(typical representation of a quadratic form).
This equation follows from the corresponding classical identity
(cf. \cite{Cl}) and Lemma \ref{Sa12}.

It follows the decomposition into two linear forms
\begin{eqnarray*}
 f \ P_{01}^2 f\ & = & (P_{01}f+(01)\sqrt{-\frac{1}{2}d_2})
                       (P_{01}f-(01)\sqrt{-\frac{1}{2}d_2}) \\
                 & =: & (xY_1-yX_1-hyY_1)(xY_2-yX_2-hyY_2).
\end{eqnarray*}
We obtain for the homogeneous coordinates
\begin{eqnarray*}
X_{1,2} & := &
-x_1B-y_1C \pm
x_1 \sqrt{-\frac{1}{2}d_2}-2hy_1B-\frac{3}{2}hx_1 A,
\\
Y_{1,2} & := &
 x_1 A+y_1B\pm y_1\sqrt{-\frac{1}{2}d_2}-\frac{1}{2}hy_1 A.
\end{eqnarray*}
For the calculations we used the central skew field extension
\[
Q_I\left[\sqrt{-\frac{1}{2} d_2}\right]\ = \
Q_I\left[\sqrt{B^2 - AC + h AB - \frac{3h^2}{4} A^2 } \right].
\]
For the corresponding projective coordinates
$Z_i=X_iY_i^{-1}+\frac{h}{2}$, $i=1,2$ we have
\[  Z_{1,2} := X_{1,2}\ Y_{1,2}^{-1} + \frac{h}{2}\ =\
      (-B \pm \sqrt{-\frac{1}{2}d_2})A^{-1}.      \]
$Z_1,Z_2$ are obviously  independent of $x_1,y_1$.

\begin{flushleft}
\bf \large 12. The solution of the cubic equation
\end{flushleft}

Let
\[  f = x^3A+x^2yB+xy^2C+y^3D           \]
be a cubic form with the discriminant $d_3$,
with the Hessian covariant $\Delta$
and with the cubic  covariant $j$,
cf. Example \ref{Ex7}.
Analogous to the classical case we have the
equation
\begin{eqnarray} \label{eq22}
    f \cdot  (P_{01}^3 f)^2\ =
    \ (P_{01}^2 f)^3\ +
    \ 3\ P_{01}^2\Delta \cdot (01)^2 P_{01}^2f\  -
    \ P_{01}^3Q \cdot (01)^3
\end{eqnarray}
(typical representation of a cubic form,
cf. \cite{Cl} and Lemma \ref{Sa12}).
With the notations $\Delta_1:=P_{01}^2 \Delta$,
$j_1:=P_{01}^3 j$,
$\xi=P_{01}^2f=P_{10}f_1$ and $\eta=(01)$ we have
\begin{eqnarray} \label{eq23}
  f\  f_1^2\ =
    \ \xi^3\ +
    \ 3\ \Delta_1\ \eta^2\ \xi\  -
    \ \ j_1 \eta^3.
\end{eqnarray}
This is a commutative cubic form in
$\xi$ and $\eta$
and we can apply the classical Cardano formula.
We make the ansatz
\begin{eqnarray} \label{eq24}
 \xi\eta^{-1} = \varepsilon^i u_1 +\varepsilon^{2i} u_2
\end{eqnarray}
with two commuting variables $u_1,u_2$ and
$\varepsilon=e^{2\pi i\over 3}$.

We insert (\ref{eq24}) into (\ref{eq23}). By Corollary
\ref{Sa6} there are no algebraic relations between $\xi ,\eta
,\Delta_1 ,j_1 $. Therefore a comparison of the coefficients
yields
\begin{eqnarray*}
u_1^3+u_2^3 & = &   j_1, \\
u_1 u_2     & = & - \Delta_1 .
\end{eqnarray*}
$u_1^3$, $u_2^3$
are the roots of the quadratic equation
\[   z^2 - j_1 z - \Delta_1^3 = 0.     \]
It follows
\[ u_{1,2}^3=\frac{j_1}{2}
\pm \sqrt{\frac{j_1^2}{4} + \Delta_1^3 }. \]
Because of
\begin{eqnarray*}
4 \Delta_1^3\ +\ j_1^2\ +\ d_3\ f_1^2\ =\ 0
\end{eqnarray*}
(cf. Remark \ref{Rem9}) we have
\[  u_{1,2}=\sqrt[3]{
\frac{1}{2} \left( j_1 \pm f_1\ \sqrt{-d_3  } \right) }. \]
Therefore we have a decomposition into three linear forms
\[  f=C( \xi - (              u_1+               u_2 )\eta )\
       ( \xi - (\varepsilon   u_1+ \varepsilon^2 u_2 )\eta )\
       ( \xi - (\varepsilon^2 u_1+ \varepsilon   u_2 )\eta )  \]
with $C\in Q^{Inv}_{I\backslash \{0\} }$.
We have
\[ \xi - (\varepsilon^i u_1+ \varepsilon^{2i} u_2 )\eta
\ =\  xY_i-yX_i-hyY_i  \]
with
\[  X_i = \gamma
            - (\varepsilon^i u_1+ \varepsilon^{2i} u_2 ) x_1, \]
\[  Y_i = \delta
            - (\varepsilon^i u_1+ \varepsilon^{2i} u_2 ) y_1 \]
($i=1,2,3$).
$\gamma ,\delta $ have the explizite form
\[ \gamma = - x_1^2 (B+hA) - x_1y_1 (2C+3hB+h^2A)
                           - y_1^2 (D+2hC+h^2B-h^3A), \]
\[ \delta =  x_1^2 A + x_1y_1 (2B-hA) + y_1^2 (C-hB+h^2A). \]
For the calculations we used the central skew field extension
\[      Q_I\left[ \ \sqrt{-d_3},\ \sqrt[3]{
\frac{1}{2} \left( j_1 \pm f_1\ \sqrt{- d_3 } \right)}
\ \right]. \]
For the projective coordinates $Z_i=X_iY_i^{-1}+\frac{h}{2}$,
$i=1,2,3$ we obtain
\[  Z_i = (\gamma
            - (\varepsilon^i u_1+ \varepsilon^{2i} u_2 ) x_1)
          (\delta
            - (\varepsilon^i u_1+ \varepsilon^{2i} u_2 ) y_1)^{-1}
+\frac{h}{2}.   \]

\begin{Rem} \label{Rem14} \rm
Without proof
we mention an explizite formula for $Z_i$
independent of $x_1,y_1$.
We extend $Q_I^\Sigma$ by the noncentral cubic root elements
$w_\pm =  \sqrt[3]{
\frac{1}{2} \left( S \pm A\ \sqrt{- d_3 } \right)}$
where $S=A^2D-3ABC+2B^3+9hAB^2-2h^2A^2B+6h^3A^3$
is the first coefficient of $j$.
We require the commutation relations
$w_+ w_- = w_- w_+$,
$w_\pm y_i = y_i w_\pm $ and
$w_\pm x_i = x_i w_\pm - \frac{1}{3}h y_i w_\pm $.
Then we have
\[Z_i = - B A^{-1}+\frac{3}{2}h +\epsilon^i A^{-1}   \sqrt[3]{
\frac{1}{2} \left( S + A\ \sqrt{- d_3 } \right)}+
    \epsilon^{2i} A^{-1} \sqrt[3]{
\frac{1}{2} \left( S - A\ \sqrt{- d_3 } \right)}\ .  \]
\end{Rem}

\begin{Rem} \label{Rem15} \rm
The quartic form has the typical representation
\[ ff_1^3= \xi^4 +3\Delta_1 \xi^2 \eta^2 + T_1 \xi\eta^3 +
\left( \frac{P f_1^2}{2}-\frac{3\Delta_1^2}{4} \right) \eta^4  \]
with the Hessian covariant $\Delta$ of degree 4, the skew
covariant $T$ of degree 6 and the fundamental invariant $P$ of
order 2 (cf. \cite{Cl}).
\end{Rem}

\begin{flushleft}
\bf \large 13. The total differential and differential modules
\end{flushleft}

In the following we introduce differentials
and differential modules.
Consider the algebra $H_I$.
We fix a subset $K\subset I$ with $0\in K$.
The elements $x_i,\ y_i$ ($i\in K$)
play the role of variables,
the elements $x_i,\ y_i$ with $i\in I\backslash K$
are constant with respect to our derivation $d$.

We define the {\it differential module}
$\Gamma_{I,K}=Q_I d_K Q_I$.
(Analogously one can consider differential modules
$H_IdH_I$.)
$\Gamma_{I,K}$
is a $Q_I$-module whose elements
are linear combinations of monomials
\[     a\ d x_i,\ \ \ \ \
       b\ d y_i,\ \ \ \ \
       d x_i \cdot c,\ \ \ \ \
       d y_i \cdot d\ \]
with $a,b,c,d \in Q_I$, $i\in K$.
(The differentials $d x_i=d_K x_i$, $d y_i=d_K y_i$ with
$i\in K$ are abstract elements).
The defining relations of the $Q_I$-module are
\[      d(\alpha a+\beta b)-\alpha da-\beta db=0,       \]
\[      d(ab)-da\cdot b-adb=0,   \]
\[           d1=0    \]
for $a,b\in Q_I$, $\alpha,\beta\in \CC$ and
\begin{eqnarray} \label{eq25}
\begin{array}{ll}
dx_j \cdot x_i -
x_i dx_j +h x_idy_j-hy_idx_j-h^2y_idy_j = 0,
\ \ \ \ \ \ \   &
dy_j \cdot y_i - y_i dy_j = 0, \\
dx_j \cdot y_i - y_i dx_j -h y_idy_j = 0, &
dy_j \cdot x_i - x_i dy_j +h y_idy_j = 0
\end{array}
\end{eqnarray}
for $i \in I$, $j\in k$.

These relations allow us to represent
every element of $\Gamma_{I,K}$ in the form
\begin{eqnarray} \label{eq26}
\sum_{i\in I} a_idx_i+b_idy_i
\end{eqnarray}
with $a_i$, $b_i\in Q_I$.

We define the {\it total differential} $d=d_K$ as a linear map
$d_K:Q_I\rightarrow \Gamma_{I,K}$ by
\begin{eqnarray*}
\begin{array}{ccl}
d_K(x_i) = dx_i, \ \ \ \ \ &
d_K(y_i) = dy_i \ \ \ \ \  &
\ \ \ \ \ \forall i\in K,  \ \ \ \       \\
d_K(x_i) = 0, \ \ \ \ \ &
d_K(y_i) = 0 \ \ \ \ \  &
\ \ \ \ \ \forall i\in I\backslash K
\end{array}
\end{eqnarray*}
and we require the Leibniz rule
\[   d_K(ab) = d_K a \cdot b + a\ d_K b  \]
for $a,b\in Q_I$.

\begin{Sa} \label{Sa16}
$d_K$ is a well defined linear map.
\end{Sa}

\noindent {\it Proof}: One has to check the consistency of the
Leibniz rule with the relations (\ref{eq4}). E.g. we have for
$i,j\in K$
\[  d(x_jy_i) = x_jdy_i +dx_j\cdot y_i
=   dy_i\cdot x_j+hdy_i\cdot y_j + y_i dx_j+h y_i dy_j
=   d(y_ix_j+hy_iy_j). \bullet        \]

\begin{Rem} \label{Rem16} \rm
We identify $K$ with a set $K'$ with $I\cap K'=0$.
Because of the similarity of the
commutation relations (\ref{eq25}) and (\ref{eq4})
we can identify $\Gamma_{I,K}$
with a linear subspace of $Q_{I\cup K'}$ according to $dx_i \cong
x_{i'}$ and $dy_i\cong y_{i'}$.
\end{Rem}

Because of the identification in Remark \ref{Rem16} we can
consider $\Gamma_{I,K}$ as a left $U_h(sl(2))$-module. The
action is given by
\begin{eqnarray*} \label{eqz2}
\begin{array}{lll}
E dx_i =  0,    \ \ \ \ \   &
F dx_i =  dy_i, \ \ \ \ \    &
H dx_i =  dx_i,   \\
E dy_i =  dx_i, \ \ \ \ \    &
F dy_i =  0,     \ \ \ \ \   &
H dy_i =  - dy_i
\end{array}
\end{eqnarray*}
for $i\in K$
and by an extension
of the Leibniz rules (\ref{eq3}) to $\Gamma_{I,K}$.

\begin{Sa} \label{Sa17}
$d_K$ commutes with the left $U_h(sl(2))$-action.
\end{Sa}

\noindent {\it Proof}:
By (\ref{eqz2}) the proposition is valid for the generators
$x_i,y_i\in H_I$.
We suppose the proposition for $a,b\in H_I$ with
$\deg(a),\deg(b)\leq n$.
Because of $d e^{hF}(a)=e^{hF}(da)$ for $a\in H_I$ we have
\[
Ed(ab) = E(da\cdot b + adb )
       = E(da)  e^{-hF} (b)  + e^{hF}(da)  E(b) +
         E(a)  e^{-hF} (db)  + e^{hF}(a)  E(db)    \]
\[     = dE(a)  e^{-hF} (b)  +d e^{hF}(a)  E(b)+
         E(a)  de^{-hF} (b)  + e^{hF}(a)  dE(b)
     = d( E(a)  e^{-hF} (b)  + e^{hF}(a)  E(b) )=dE(ab).
\]
Similar we obtain $dF=Fd$, $dH=Hd$.
The proposition follows.$\bullet$

\begin{De} \label{De4} \rm
We call $\gamma \in \Gamma_{I,K}$ a {\it
differential invariant}, if $E\gamma =F\gamma =H\gamma =0$.
We denote the submodule of differential invariants
by $\Gamma_{I,K}^{Inv}$.
\end{De}

By Proposition \ref{Sa17} we have
\[ d_K(Q_I^{Inv})\subseteq \Gamma_{I,K}^{Inv}. \]
For $i\in K$, $j\in I\backslash K$
we obtain the differential invariant
\[  d_K( (ij) ) =  x_idy_j-y_idx_j-hy_idy_j.  \]
We will use the notation $(idj)=x_idy_j-y_idx_j-hy_idy_j$.

We need the following simple Lemma about the
structure of differential modules.
\begin{Le} \label{Sa18}
Let $K,K'\subset I$ and $K\cap K'=\emptyset$.
Then we have
$\Gamma_{I,K\cup K'}=\Gamma_{I,K}\oplus\Gamma_{I,K'}$
and $d_{K\cup K'}=d_K + d_{K'}$.
\end{Le}
The proof follows from the Leibniz rule.

\begin{Rem} \label{Rem17} \rm
Differential modules $Sd_KS$
for central field extensions
$S=Q_I[\sigma_1,....,\sigma_k]$ are defined analogous
to $Q_IdQ_I$.
Additional requirements are the defining relations of the $\sigma_i$
and the conditions
$\sigma_i dx-dx\cdot \sigma_i = x d\sigma_i-d\sigma_i\cdot x=0$
for $x\in S$.
\end{Rem}

We can extend the polarization operators $P_{ij}$ and
$P_{(x_i,y_i),(X,Y)}$ to $H_IdH_I$. The extensions are given if we
set $P_{ij}(da):=dP_{ij}(a)$,
$P_{(x_i,y_i),(X,Y)}(da):=dP_{(x_i,y_i),(X,Y)}(a)$ and require
that the Leibniz rule (\ref{eq14}) is valid for $a,b\in H_IdH_I$.
By Remark \ref{Rem16} the extended operators are well
defined.

\begin{flushleft}
{\bf \large 14. Polynomials }
\end{flushleft}

\noindent
In this section we consider polynomials
associated to $n$-forms.
\begin{Le} \label{Sa19}
The elements $y_i$, $(ij)$ and $z_i-z_j$ generate a commutative
subfield of $Q_I$.
\end{Le}
This follows immediately from $y_iy_j=y_jy_i$,
$z_i-z_j=(ij)y_i^{-1}y_j^{-1}$ and
Proposition \ref{Sa3}.

We introduce the bracket
\[  [ij):=y_i^{-1}(ij)=(ij)y_i^{-1}=(z_i-\frac{h}{2})y_j-x_j \]
with the properties $[ij)[kl) = [kl)[ij)$
and $[ij)y_k  = y_k [ij)$
for all $i,j,k,l$. Let $f=(01)(02)...(0n)$ be a $n$-form. We
introduce the the polynomial $f_z$ of $f$ by
\[ f_z:=y^{-n} f = [01)[02)...[0n)
\ =\ \prod_{i=1}^n \left( (z-\frac{h}{2})y_i-x_i \right).  \]
By polarization we obtain
\[ f_Z: = Y^{-n} P^n_{(x,y),(X,Y)} f = [Z1)[Z2)...[Zn)
\ =\ \prod_{i=1}^n \left( (Z-\frac{h}{2})y_i-x_i \right) \] with
$Z:=XY^{-1}+\frac{h}{2}$ and $[Zi):=(Z-\frac{h}{2})y_i-x_i$.

We consider the differential module
$\Gamma_{I,K}=Q_Id_KQ_I$.
We have
\[  d_i[ij)=y_j dz_i\ \ \ \ \ \ \ \ \ \  {\rm if}
\ \ \ \ \ \ \ i\in K,j\notin K     \]
and
\[  d_j[ij)=(z_i-\frac{h}{2})dy_j-dx_j=:[0dj)
\ \ \ \ \ \ \ \ \ \
{\rm if}
\ \ \ \ \ \ \ i\notin K,j\in K.   \]
Furthermore we have
\[  dz = dx\cdot y^{-1}+xdy^{-1} = dx\cdot y^{-1}-xy^{-2}dy =
(dx\cdot y-dy\cdot x - hdy\cdot y)y^{-2}= - (0d0)y^{-2}.   \]
We generalize Lemma \ref{Sa19}.
\begin{Le} \label{Sa20}
The elements $y_i,z_i-z_j,(ij),[ij),dy_i,dz_i,(idj),[idj)$ commute
in $\Gamma_{I,K}$ .
\end{Le}
We compute the differentials of $f_z$.
Let $K=\{ 0,1,2,...,n \}$ and
$d=d_0+\delta =d_0+d_{\{1,2,...,n \}}$.
\begin{Le} \label{Sa21}
We have
\[ d_0f_z= d_0 ([01)[02)...[0n)) =
\left( \sum_{i=1}^n y_i [01)...[0,i-1)[0,i+1)...[0n) \right)dz \]
and
\[ \delta f_z= \delta ([01)[02)...[0n)) =
\sum_{i=1}^n [01)...[0,i-1)[0di)[0,i+1)...[0n). \]
\end{Le}

\begin{flushleft}
\bf \large 15. Elliptic and Hyperelliptic Differentials
\end{flushleft}

We extend $Q_I$ by the central element $w$ with
\[  w^2 = (01)(02)(03)(04).      \]
We obtain the central skew field extension $Q_I[w]$.
We denote $w$ by
\[  w = \sqrt{(01)(02)(03)(04)}.  \]
We can suppose that $w$ is invariant.
\begin{De} \label{De5} \rm
We say that the differential invariant
\begin{eqnarray} \label{eq27}
de_{x,y}:= w^{-1} y^{2} dz =
\sqrt{(01)(02)(03)(04)}^{-1}\ \ (ydx-xdy+hydy)
\end{eqnarray}
of the differential module $Q_I[w]d_K Q_I[w]$
with $0\in K$ and $1,2,3,4\notin K$
is an {\it elliptic differential of the first kind}.
\end{De}

\begin{Rem} \label{Rem18}\rm
In the classical case we have
{\small
\[  de_{x,y}
=\frac{ydx-xdy}{
\sqrt{(xy_1-yx_1)(xy_2-yx_2)(xy_2-yx_2)(xy_2-yx_2)}}
=\frac{\frac{1}{\sqrt{y_1y_2y_3y_4}}\ \ \  dz
}{\sqrt{(z-z_1)(z-z_2)(z-z_3)(z-z_4)}}. \]}
\end{Rem}

More generally, for $g\in \NN$, $g\geq 1$
we extend $Q_I$ by the central
and invariant element $W$
with
\[  W^2 = (0i_1) \cdots (0i_{2g+2})      \]
where the indices $i_j$ are different.
We obtain the central skew field extension $Q_I[W]$.
We use the notation
\[  W = \sqrt{(0i_1) \cdots (0i_{2g+2})}. \]
Furthermore we consider the $(g-1)$-form
\[  U = (0j_1) \cdots (0j_{g-1}).  \]
\begin{De} \label{De6} \rm
We say that the differential invariant
\begin{eqnarray} \label{eq28}
\begin{array}{c}
 dh_{x,y}:=  U W^{-1} y^{2} dz   \\ =
(0j_1) \cdots (0j_{g-1})\ \
\sqrt{(0i_1) \cdots (0,i_{2g+2})}^{-1}\ \ (ydx-xdy-hydy)
\end{array}
\end{eqnarray}
of the differential module $Q_I[W]d_K Q_I[W]$
with $0\in K$ and $i_1,...,i_{2g+2},j_1,...,j_{g-1}\notin K$
is a {\it hyperelliptic differential of the first kind}.
\end{De}

\begin{Rem} \label{Rem19}\rm
In the classical case we have
{\small \[  dh_{x,y}
=\frac{(xy_{j_1}-yx_{j_1}) \cdots (xy_{j_{g-1}}-yx_{j_{g-1}})
\ \ \ (ydx-xdy)}{
\sqrt{(xy_{i_1}-yx_{i_1})(xy_{i_2}-yx_{i_2})
    .............. (xy_{i_{2g+2}}-yx_{i_{2g+2}})}} \]
\[
=\frac{y_{j_1}...y_{j_{g-1}}}{ \sqrt{y_{i_1} \cdots y_{i_{2g+2}}}} \ \ \
 \frac{
 (z-z_{j_1})\cdots (z-z_{j_{g-1}}) \ \ \ dz}{
 \sqrt{(z-z_{i_1})(z-z_{i_2}) ............. (z-z_{i_{2g+2}})}}.
\]}
$dh_{x,y}$ is connected with an algebraic curve of genus $g$. In
particular for $g=1$ we obtain elliptic differentials. The $g$
classical hyperelliptic differentials of the first kind $\frac{z^i
dz}{\sqrt{p(z)}}$, $i=0,1,2,...,g-1$ correspond to $dh_{x,y}$ with
$U=(0j_1)^i (0j_2)^{g-1-i}$,
where $(x_{j_1},y_{j_1})$ is the zero point
 and  $(x_{j_2},y_{j_2})$ is the point
at infinity.
\end{Rem}

We can extend the polarization process $P_{(x,y),(X,Y)}$
to extensions of fields and of differential modules.
We change the arguments of elliptic
and hyperelliptic differentials  by the definition
\[ de_{X,Y}:=
(\sqrt{P_{(x,y),(X,Y)}^4  w^2})^{-1} P_{(x,y),(X,Y)}^2 (d0,0) \]
and
\[ dh_{X,Y}:=
(\sqrt{P_{(x,y),(X,Y)}^{2g+2} W^2})^{-1}
P_{(x,y),(X,Y)}^{g+1} (U\ (d0,0)).          \]
for arbitrary noncommutative points $(X,Y)$
with $Z=XY^{-1}+\frac{h}{2}$.

\begin{Rem} \label{Rem20} \rm
We obtain $de_{X,Y}$ and $dh_{X,Y}$
by replacing $(x,y)$ by $(X,Y)$
in $de_{x,y}$ and $dh_{x,y}$, respectively.
\end{Rem}

\begin{flushleft}
\bf \large 16. The Addition theorem and Abel's Theorem
\end{flushleft}

Let $\Gamma_{I,K}$ be a differential module
with $1,2,3,4 \notin K$ and $0,5 \in K $.

For the above elliptic differential (\ref{eq27})
we consider the cubic form
\[  r =(01)(02)(03)-(04)(05)^2.     \]
According to Theorem \ref{Sa14}
we have the decomposition into commuting factors
\[  r = C (xY_1-yX_1-hyY_1)(xY_2-yX_2-hyY_2)(xY_3-yX_3-hyY_3) \]
with $X_i,Y_i$ in a splitting field
$Q_{I\backslash \{0\} }^\Sigma  $ and
$C\in Q^{Inv}_{I\backslash \{0\}}$.
We have
$ r_{X_i,Y_i}:=P^{3}_{(x,y),(X_i,Y_i)} r = 0$
for $i=1,2,3$.

Furthermore we extend
$Q_{I}^\Sigma $
by the central root elements
\[ w^{(i)} := \sqrt{ P^4_{(x,y),(X,Y)} w^2 } =
\sqrt{(X_i1)(X_i2)(X_i3)(X_i4)}, \ \ \ \ \ i=1,2,3  \]
with  $(X_ij):= X_iy_j-Y_ix_j-hY_iy_j $
to the skew field
$Q_{I}^{\Sigma ,w} $.

Now we consider the elliptic differentials
\[    de_{X_i,Y_i}:= (w^{(i)})^{-1} (dX_i,X_i)      \]
of the differential module
$\Gamma_{I,K}:=
Q_{I}^{\Sigma ,w}
d_{K} Q_{I}^{\Sigma ,w} $.

We formulate the addition theorem.
\begin{Th} \label{Sa22}
Let $(X_i,Y_i)$ be the homogeneous coordinates
of the linear factors of $r=(01)(02)(03)-(04)(05)^2$.
Then we have
\[  \epsilon_1 de_{X_1,Y_1}+
    \epsilon_2 de_{X_2,Y_2}+
    \epsilon_3 de_{X_3,Y_3} = 0   \]
with certain $\epsilon_i =\pm 1 $.
\end{Th}

\begin{Rem} \label{Rem21} \rm
In the classical case the Theorem reduces to the differential
equation
{\tiny
\[ \frac{\ \ \  dZ_1
}{\sqrt{(Z_1-z_1)(Z_1-z_2)(Z_1-z_3)(Z_1-z_4)}} + \frac{\ \ \  dZ_2
}{\sqrt{(Z_2-z_1)(Z_2-z_2)(Z_2-z_3)(Z_2-z_4)}} + \frac{\ \ \  dZ_3
}{\sqrt{(Z_3-z_1)(Z_3-z_2)(Z_3-z_3)(Z_3-z_4)}}=0
\]}
if $Z_1,Z_2,Z_3$ vary continuously as the zeros
of the equations
\[
(z-z_1)(z-z_2)(z-z_3)-(z-z_4)(az+b)^2=0  \]
with arbitrary $a,b\in \CC$.
We can consider
$Z_1,Z_2,Z_3$ as the coordinates of
the intersection points of the lines $q(z,u,a,b)=u-az-b=0$ and the
elliptic curve \\
$p(z,u)=u^2(z-z_4)-(z-z_1)(z-z_2)(z-z_3)=0$
or
as the intersection points of the variable curve
$q(z,u,a,b)=(az+b)u-(z-z_1)(z-z_2)(z-z_3)=0$
and the
elliptic curve
$p(z,u)=u^2-(z-z_1)(z-z_2)(z-z_3)(z-z_4)=0$.
If we choose two different points
with coordinates $Z_1,Z_2$
and $\epsilon_1,\epsilon_2$ then
the parameters $a,b$ are fixed
and $Z_3$, $\epsilon_3$
are uniquely determined (classical addition theorem).
\end{Rem}

For the consideration of hyperelliptic differentials
we suppose
$i_1,i_2,\cdots ,i_{2g+2},j_1,\cdots ,j_{g-1} \notin K$ and
$0,k_1,\cdots ,k_p,l_1,\cdots ,l_q \in K$.
We require that the indices $i_j,k_j,l_j$ are different
(cf. below for $k_j,l_j$).

For the above hyperelliptic differential
we fix a number $s=0,1,...,2g+2$
and two forms
$P=(0k_1)....(0k_p)$,
$Q=(0l_1)....(0l_q)$
of degree $p$ and $q$, respectively,
with
\[  p - q = g + 1 - s.   \]
It follows
\[  s+2p=2g+2-s+2q.  \]

\begin{Rem} \rm The equation $p-q=g+1-s$ is a technical condition
for this invariant theoretical consideration
in order to secure the homogeneity of $r$ (cf. below).
Furthermore one can replace $P$ and $Q$ by arbitrary
$p$- and $q$-forms, respectively.
\end{Rem}

Furthermore we consider the decomposition
\[   W^2 = AB \]
with $A=(0i_1)....(0i_s)$ and $B=(0,i_{s+1})...(0,i_{2g+2})$.
Therefore
\[  r:=(0i_1)....(0i_s)P^2-(0,i_{s+1})...(0,i_{2g+2})Q^2
=  A P^2-B Q^2       \]
is the difference of two
$k$-forms with $k:=s+2p=2g+2-s+2q=p+q+g+1$.

According to Theorem \ref{Sa14}
we have the decomposition
\[  r = C(0X_1)\cdots (0X_k)     \]
into commuting factors $(0X_i):= xY_i-yX_i-hyY_i$
with $X_i,Y_i$ in a splitting field
$Q_{I\backslash \{0\} }^\Sigma  $ and
$C\in Q^{Inv}_{I\backslash \{0\} }$.
We have
$r_{X_i,Y_i}=P_{(x,y),(X_i,Y_i)}^k r =0$ for
$i=1,\cdots k$.

Furthermore we extend
$Q_{I}^\Sigma $ by the  root elements
\[  W^{(i)} :=  \sqrt{ P_{(x,y),(X_i,Y_i)}^{2g+2} W^2 }=
               \sqrt{(X_ii_1)\cdots (X_i i_{2g+2})}  \]
to $Q_{I}^{\Sigma ,W}$.

Then we consider the differential module
$\Gamma_{I,K}=
Q_{I}^{\Sigma ,W} d_K
Q_{I}^{\Sigma ,W}$ and
the hyperelliptic differentials
\[ dh_{X_i,Y_i} = {W^{(i)}}^{-1}
   P_{(x,y),(X_i,Y_i)}^{g-1}U\  (dX_i,X_i) \]
with $i=1,2,...,k$.
We formulate Abel's Theorem:
\begin{Th} \label{Sa23}
Let $(X_i,Y_i)$, $i=1,2,...,k$
be the homogeneous coordinates
of the linear factors of
$r=(0i_1)....(0i_s)P^2-(0,i_{s+1})...(0,i_{2g+2})Q^2$.
Then we have
\[  \epsilon_1 dh_{X_1,Y_1}+.....+
    \epsilon_k dh_{X_k,Y_k} = 0      \]
with certain $\epsilon_i =\pm 1 $.
\end{Th}

\begin{Rem} \label{Rem22} \rm
In the classical case we have
\[    \sum_{i=1}^k
\frac{U_z(Z_i)dZ_i}{\sqrt{(Z_i-z_1)(Z_i-z_2)...(Z_i-z_{2g+2})}}
=0
\]
where $Z_1,...,Z_k$ vary continuously as coordinates of the
intersection points of the curves $c(z,u,a_0,...,a_q,b_0,...,b_p)=
(a_0+a_1z+...+a_qz^q)u-(b_0+b_1z+...+b_pz^p)(z-z_1)...(z-z_s)=0$
and the hyperelliptic curve $p(z,u)=u^2 -
(z-z_{1})(z-z_{2})...(z-z_{2g+2})=0$.
A point of $p(u,z)=0$ is given by a pair $(Z,\epsilon)$
with $\epsilon=\pm 1$.
We choose different
points $(Z_1,\epsilon_1),...,(Z_{p+q+1},\epsilon_{p+q+1})$.
Then we obtain a
system of the $p+q+1$ linear equations
\[    c(Z_i, \epsilon_i \sqrt{( Z_i-z_1)...(Z_i-z_{2g+2})},
   a_0,...,a_q,b_0,...,b_p) = 0 \]
for the p+q+2 coefficients $a_0,...,b_p$.
In the generic case the rank  is $p+q+1$
and the coefficients are determined
up to a common factor.
In this non-degenerate case the remaining $g$ points
$(Z_i,\epsilon_i)$ are uniquely determined.
\end{Rem}

\begin{flushleft}
\bf \large 17. Proof of Theorem \ref{Sa23}
\end{flushleft}

We prove Theorem \ref{Sa23}.
We obtain Theorem \ref{Sa22}
as the special case $g=1,s=3,p=0$,
(i.e. $q=1,k=3,W=w,A=(01)(02)(03),B=(04),U=1,P=1,Q=(05)$).

Using invariant theory,
the the proof is analogous to the classical case
(cf. \cite{Ab,We}).
We consider the $k$-form
$r=A P^2-B Q^2$
with the decomposition \\
$r=C (xY_1-yX_1-hyY_1)...(xY_k-yX_k-hyY_k)$
with $X_i,Y_i$ in a splitting field
$Q_{I\backslash \{0\} }^\Sigma  $ and
$C\in Q^{Inv}_{I\backslash \{0\} }$ and with
$r_{X_i,Y_i}=P_{(x,y),(X_i,Y_i)}^k r =0$
for $i=1,...,k$.

In order to apply $d=d_K$ we introduce $r_z$.
We have
\[ r_z = y^{-k} r = C [0 X_1) ... [0 X_k)
\ =\ \prod_{i=1}^k \left( (z-\frac{h}{2})Y_i-X_i \right).  \]
We consider the differential
$d=d_0+\delta := d_0 + {d}_{\{k_1,...k_p,l_1,...,l_q \} }$
(cf. Lemma \ref{Sa18}).
We obtain
\[  dr_z = d_0 r_z + \delta r_z =
d_0 r_z +2 A_z P_z \delta P_z - 2 B_z Q_z \delta Q_z  \]
with
\[   d_{0} r_{z} = \sum_{i=1}^k
 C Y_i[0 X_1)...[0 X_{i-1})[0 X_{i+1})
 ...[0 X_k) dz ,       \]
\[    \delta P_{z}=
      \sum_{i=1}^p [0 k_1)...[0 k_{i-1})[0 dk_i)
           [0 k_{i+1})...[0 k_p)    \]
and
\[ \delta Q_{z}=
   \sum_{i=1}^q [0 l_1)...[0 l_{i-1})[0 dl_i)
   [0 l_{i+1})...[0 l_q),      \]
cf. Lemma \ref{Sa21}. We multiply both sides with
$y^{k+1}=y^{p+q+g+2}$. Using $y \circ \delta = \delta \circ y $
we obtain
\begin{eqnarray} \label{eq30}
y^{k+1}  dr_{z}=
y^{k+1}  d_{0} r_{z} +
2 y A P \delta P -
2 y B Q \delta Q
\end{eqnarray}
with
\[y^{k+1}  d_{0} r_{z} =
\sum_{i=1}^k
 C Y_i(0 X_1)...(0 X_{i-1})(0 X_{i+1})
  ... (0 X_k) (d0,0),       \]
\[    \delta P =
      \sum_{i=1}^p (0 k_1)...(0 k_{i-1})(0 dk_i)
           (0 k_{i+1})...(0 k_p)    \]
and
\[    \delta Q =
          \sum_{i=1}^q (0 l_1)...(0 l_{i-1})(0
          dl_i)
           (0 l_{i+1})...(0 l_q).      \]

We apply $P^{k+1}_{(x,y),(X_i,Y_i)}$ to
(\ref{eq30}).
Because of $r_{X_i,Y_i}=0$ we obtain for the left side
\[   P^{k+1}_{(x,y),(X_i,Y_i)}
     (y^{k+1}
     d(y^{-k}r ))=
     P^{k+1}_{(x,y),(X_i,Y_i)}
     ( - kdy\cdot r + y dr) =
       - kdY_i\cdot r_{X_i,Y_i} + Y_i dr_{X_i,Y_i} = 0.\]
Therefore we have for the right side
\begin{eqnarray} \label{ll}
\begin{array}{c}
0  =
 C Y_i(X_i X_1)...(X_i X_{i-1})(X_i X_{i+1})
  ... (X_i X_k)  (dX_i,X_i) + \\
2Y_i A_{X_i,Y_i } P_{X_i,Y_i}
                  \delta P_{X_i,Y_i } -
2Y_i B_{X_i,Y_i } Q_{X_i,Y_i}
                  \delta Q_{X_i,Y_i }.
\end{array}
\end{eqnarray}

We have for $i=1,2,...,k$ the identity
\[ r_{X_i,Y_i} =  P_{(x,y),(X_i,Y_i)}^k r =0. \]
It follows
\begin{eqnarray} \label{eqz3}
  A_{X_i,Y_i}   P_{X_i,Y_i}^2 - B_{X_i,Y_i}    Q_{X_i,Y_i}^2 = 0,
\end{eqnarray}
\[  A_{X_i,Y_i}^2 P_{X_i,Y_i}^2 - (W^{(i)})^2  Q_{X_i,Y_i}^2 = 0,  \]
\[ (A_{X_i,Y_i}   P_{X_i,Y_i} + W^{(i)}    Q_{X_i,Y_i})
   (A_{X_i,Y_i}   P_{X_i,Y_i} - W^{(i)}    Q_{X_i,Y_i})=0. \]
Because the skew field $Q_{I }^{\Sigma ,W} $
has no zero divisors it follows
\begin{eqnarray} \label{eq31}
A_{X_i,Y_i}     P_{X_i,Y_i} = \epsilon_i W^{(i)}  Q_{X_i,Y_i}
\end{eqnarray}
with $\epsilon_i = 1$ or $-1$.
From (\ref{eqz3}) and (\ref{eq31}) it follows
\begin{eqnarray} \label{eq32}
B_{X_i,Y_i}     Q_{X_i,Y_i} = \epsilon_i W^{(i)}  P_{X_i,Y_i}.
\end{eqnarray}

We multiply (\ref{ll}) with $Y_i^{-1}$
and we insert the relations (\ref{eq31}) and (\ref{eq32}).
\begin{eqnarray*}
 C (X_i X_1)...(X_i X_{i-1})(X_i X_{i+1})
  ... (X_i X_k)  (dX_i,X_i)
=  -
2 \epsilon_i W^{(i)} Q_{X_i,Y_i}
                  \delta P_{X_i,Y_i }
+
2 \epsilon_i W^{(i)} P_{X_i,Y_i}
                  \delta Q_{X_i,Y_i }.
\end{eqnarray*}

Because $i_j,k_j,l_j$ are different indices,
$r$ has $k$ different zeros.
Therefore we have \\ $(X_iX_j)\neq 0$
$\forall i,j\in \{1,...,k\}$.
It follows
\[ \epsilon_i {W^{(i)}}^{-1}Y_i^{2}dZ_i =
-2C^{-1} ((X_iX_1)...(X_iX_{i-1})(X_iX_{i+1})...(X_iX_k))^{-1}
  ( Q_{X_i,Y_i} \delta P_{X_{i},Y_{i}}
 -  P_{X_i,Y_i} \delta Q_{X_{i},Y_{i}}).   \]
We multiply both sides with the polars
$U_{X_i,Y_i} = P^{g-1}_{(x,y),(X_i,Y_i)}U$
of an arbitrary $(g-1)$-form $U$.
\[ \epsilon_i dh_{X_i,Y_i}=
\epsilon_i U_{X_i,Y_i} {W^{(i)}}^{-1}Y_i^{2}dZ_i = \]
\[  -2 C^{-1} ((X_iX_1)...(X_iX_{i-1})(X_iX_{i+1})...(X_iX_k))^{-1}
U_{X_i,Y_i}  ( Q^{(i)} \delta P_{X_{i},Y_{i}}
        -  P^{(i)} \delta Q_{X_{i},Y_{i}}).   \]
We have to show that
\[   \sum_{i=1}^k \epsilon_i dh_{X_i,Y_i} = \]
\[ -2 C^{-1}      \sum_{i=1}^k
 ((X_iX_1)...(X_iX_{i-1})(X_iX_{i+1})...(X_iX_k))^{-1}
U_{X_i,Y_i}  ( Q^{(i)}
            \delta P_{X_{i},Y_{i}}
 -  P^{(i)} \delta Q_{X_{i},Y_{i}}) = 0, \]
where
\[  \delta P_{X_{i},Y_{i}} =
 \sum_{j=1}^p (X_i k_1)...(X_i k_{j-1})(X_i dk_j)
           (X_i k_{j+1})...(X_i k_p)    \]
and
\[  \delta Q_{X_{i},Y_{i}} =
 \sum_{j=1}^q (X_i l_1)...(X_i l_{j-1})(X_i dl_j)
           (X_i l_{j+1})...(X_i l_q).  \]
The proposition follows from the following Lemma.

\begin{Le} \label{Sa24}
Let $f=(0X_1)...(0X_k)$ be a $k$-form with different zeros
(i.e. $Z_i \neq Z_j$ for $i\neq j$)
and let $g$ be a $(k-2)$-form. Then we have
\[  \sum_{i=1}^{k}
\frac{ P_{(x,y),(X_i,Y_i)}^{k-2}
g }{ (X_iX_1)...(X_iX_{i-1})(X_iX_{i+1})...(X_iX_k)}\
 =\  0\ .\]
\end{Le}

\noindent {\it Proof}:
The form $g$ has a decomposition
$g=C(0{X'}_1)...(0{X'}_{k-2})$ with bracket symbols.
At first we consider
the classical case $h=0$.
Let $f_z(z)=y^{-k}f(x,y)$ and $g_z(z)=y^{-k+2}g(x,y)$.
We have the partial fraction expansion
\[    \frac{-zg_z(z)}{f_z(z)}\ =\
       \sum_{i=1}^k
       \frac{-Z_ig_z(Z_i)}{{f'}_z(Z_i)(z-Z_i)}\ \  .\]
For $z=0$ we obtain
\[  0 \ =\
       \sum_{i=1}^k
       \frac{g_z(Z_i)}{{f'}_z(Z_i)}\ =\
       \sum_{i=1}^k
\frac{g_z(Z_i)}{Y_i[Z_iX_1)...[Z_iX_{i-1})
[Z_iX_{i+1})...[Z_iX_k)} \]
\[   =  \sum_{i=1}^k
\frac{P_{(x,y),(X_i,Y_i)}^{k-2} g }{(X_iX_1)...(X_iX_{i-1})
(X_iX_{i+1})...(X_iX_k)}\  .  \]
In the special case
$X_i=x_{\alpha_i}$,
$Y_i=y_{\alpha_i}$,
${X'}_i=x_{\beta_i}$,
${Y'}_i=y_{\beta_i}$
with different $\alpha_i,\beta_i\in I$ the Lemma follows
from Lemma \ref{Sa12}.
In the general case we apply the field homomorphism with
$x_{\alpha_i}\rightarrow X_i  $,
$y_{\alpha_i}\rightarrow Y_i  $,
$x_{\beta_i}\rightarrow  {X'}_i $,
$y_{\beta_i}\rightarrow  {Y'}_i $,
$x_k \rightarrow x_k$ and
$y_k \rightarrow y_k$ for $k\neq \alpha_i,\beta_i$.
$\bullet$

We identify the differentials $dx_{k_i}$, $dy_{k_i}$, $dx_{l_i}$,
$dy_{l_i}$ with the coordinates $x_{m_i}$, $y_{m_i}$, $x_{n_i}$,
$y_{n_i}$, respectively, with $m_i,n_i\notin I$, cf. Remark
\ref{Rem16}. Therefore we can consider
 $U  ( Q  \delta P  -  P  \delta Q )$
as a $(k-2)$-form $g$.
Furthermore let $f$ be the $k$-from $r$.
The Theorem follows from Lemma \ref{Sa24}. $\bullet$


\begin{thebibliography}{9}

 \bibitem[1]{Ab}  Abel, N.H.: Remarques sur quelques proprietes
                  generales d'une certain sorte de fonctions
                  transcendantes,
                  {\it J. Reine Angew. Math.}
                  {\bf 3} (1828) 313-323, or in:
                  {\it Oeuvres}, ed. by Sylow and Lie,
                  Christiania 1881,
                  T.1, 444-456.

 \bibitem[2]{CG}  Clebsch, A. and Gordan, P.:
                  {\it Theorie der Abelschen Funktionen},
                  Teubner-Verlag, Leipzig 1866.

 \bibitem[3]{Cl}  Clebsch, A.:
                  {\it Theorie der bin\"aren algebraischen Formen},
                  Teubner, Leipzig, 1872.

 \bibitem[4]{Co}  Cohn, P.M., Skew fields, Cambridge University
                  Press, 1995.

 \bibitem[5]{Cu}  Curtis, C.W.: A note on noncommutative polynomials,
                  {\it Proc. Amer. Math. Soc.} {\bf 3} (1952), 965-969.

 \bibitem[6]{Ha}  Harris, J.: {\it Algebraic geometry: a first course},
                  Springer-Verlag, Berlin 1992.

 \bibitem[7]{He}  Hermite, Ch.: Sur la theorie des fonctions
                  homogenes a deux indeterminees,
                  {\it J. Reine Angew. Math.}
                  {\bf 52} (1866), 18-38.

 \bibitem[8]{Ka}  Kassel, C.:  {\it Quantum Groups},
                  Springer Verlag, New York, 1995.

 \bibitem[9]{Kl}  Klimek, S. and Lesniewski, A.:
                  A Two-Parameter Quantum Deformation of the Unit Disc,
                  {\it J. Func. Anal.} {\bf 115} (1993), 1-23.

 \bibitem[10]{Ko}  Kohn, G.: Zur Theorie der associirten Formen,
                  {\it Wien. Ber. C.} 1891, 865-893.

 \bibitem[11]{Kon} Kontsevich, M.:
                  Deformation quantisation of algebraic varieties,
                  {\it Preprint} AG/0106006.

 \bibitem[12]{Le} Leitenberger, F.:
                  A Quantum Deformation of Invariants of Higher
                  Binary Forms,\\
                  {\it J. Algebra}
                  {\bf 222} (1999), 82-128.

 \bibitem[13]{Ma} Majid, S.:
                  {\it Foundations of Quantum Group Theory},
                  Cambridge University Press, New York, 1995.

 \bibitem[14]{Man} Manin, Y.: Quantized Theta functions,
                  {\it Progress of Theor. Phys. Supplement}, 102 (1990)
                  219-228.

 \bibitem[15]{Mc} McConnel, J.C. and Robson, J.C.:
                  {\it Noncommutative Noetherian Rings},
                  John Wiley $\&$ Sons, Chichester, 1987.

 \bibitem[16]{Mt} Matthiessen, L.:
                  {\it Grundz\"uge der antiken und modernen Algebra
                  der litteralen Glei\-chun\-gen},
                  Teubner, Leipzig, 1878.

 \bibitem[17]{Mu} Mumford, D.: Tata lectures on theta II,
                  Birkh\"{a}user, Boston 1994.

 \bibitem[18]{Or} Ore, O.: Linear equations in non-commutative fields,
                  {\it Ann. Math.} {\bf 32} (1933), 480-508.

 \bibitem[19]{Ol} Olver, P.: {\it Classical Invariant theory},
                  Cambridge University Press, Cambridge, 1999.

 \bibitem[20]{Vb} Stafford, J.T. and Van den Bergh M.: Noncommutative
                  curves and noncommutative surfaces,
                  {\it Bull. AMS } {\bf 38}(2001)171-216.

 \bibitem[21]{St} Sturmfels, B.:
                  {\it Algorithms in Invariant Theory},
                  Springer-Verlag, Wien, 1993.

 \bibitem[22]{We} Weber, H.:
                  {\it Elliptische Funktionen und algebraische Zahlen},
                  Vieweg, Braunschweig, 1908.

\end{thebibliography}
\end{document}